\documentclass[ejs]{imsart}
\usepackage{amsmath,amstext}  
\usepackage{amsfonts,amssymb} 
\usepackage[utf8]{inputenc}
\usepackage[T1]{fontenc}
\usepackage[english]{babel}
\usepackage{graphics}
\usepackage{epsfig}
\usepackage{epsf}
\usepackage{float}
\usepackage{color}
\usepackage{pifont}
\usepackage{bbm}
\usepackage{hyperref}
\startlocaldefs
\numberwithin{equation}{section}

\newcommand{\ind}[1]{\mathbbm{1}_{\{#1\}}}
\newcommand{\indnex}{\mathbbm{1}_{\overline{\cE}}}

\newcommand{\dE}{\mathbb{E}}
\newcommand{\dF}{\mathbb{F}}
\newcommand{\dG}{\mathbb{G}}
\newcommand{\dN}{\mathbb{N}}
\newcommand{\dP}{\mathbb{P}}
\newcommand{\dR}{\mathbb{R}}
\newcommand{\dT}{\mathbb{T}}

\newcommand{\wh}{\widehat}
\newcommand{\rI}{\bs{\mathrm{I}}}
\newcommand{\rJ}{\bs{\mathrm{J}}}
\newcommand{\veps}{\varepsilon}

\newcommand{\cA}{\mathcal{A}}
\newcommand{\cB}{\mathcal{B}}
\newcommand{\cE}{\mathcal{E}}
\newcommand{\cF}{\mathcal{F}}
\newcommand{\cG}{\mathcal{G}}
\newcommand{\cH}{\mathcal{H}}
\newcommand{\cK}{\mathcal{K}}
\newcommand{\cL}{\mathcal{L}}
\newcommand{\cN}{\mathcal{N}}
\newcommand{\cO}{\mathcal{O}}
\newcommand{\cR}{\mathcal{R}}
\newcommand{\cT}{\mathcal{T}}
\newcommand{\cV}{\mathcal{V}}
\newcommand{\cW}{\mathcal{W}}
\newcommand{\cZ}{\mathcal{Z}}

\newcommand{\bs}[1]{\boldsymbol{#1}}
\newcommand{\reff}[1]{(\ref{#1})}

\def\build#1_#2^#3{\mathrel{\mathop{\kern 0pt#1}\limits_{#2}^{#3}}}
\def\liml{\build{\longrightarrow}_{}^{\cL}}
\endlocaldefs

\begin{document}
\newtheorem{Remark}{Remark}[section]
\newtheorem{Proposition}[Remark]{Proposition}
\newtheorem{Theorem}[Remark]{Theorem}
\newtheorem{Lemma}[Remark]{Lemma}
\newtheorem{Corollary}[Remark]{Corollary}
\begin{frontmatter}
\title{Parameters estimation for asymmetric bifurcating autoregressive processes with missing data}
\runtitle{Estimation for missing data BAR}

\author{\fnms{Beno\^\i te} \snm{de~Saporta}\corref{}\ead[label=e2]{saporta@math.u-bordeaux1.fr}}
\address{Universit\'e de Bordeaux, GREThA CNRS UMR~5113, IMB CNRS UMR~5251\\ and INRIA~Bordeaux~Sud~Ouest team CQFD, France\\ \printead{e2}}
\and
\author{\fnms{Anne} \snm{G\'egout-Petit}\ead[label=e3]{anne.petit@u-bordeaux2.fr}}
\address{Universit\'e de Bordeaux, IMB, CNRS UMR~525\\
 and INRIA~Bordeaux~Sud~Ouest team CQFD, France\\ \printead{e3}}
\and
\author{\fnms{Laurence} \snm{Marsalle}\ead[label=e1]{laurence.marsalle@univ-lille1.fr}}
\address{Universit\'e de Lille 1, Laboratoire Paul Painlev\'e,
CNRS UMR~8524, France\\ \printead{e1}}
 \runauthor{ B. de~Saporta, A. G\'egout-Petit, L. Marsalle}

\begin{abstract}\quad
We estimate the unknown parameters of an asymmetric bifurcating autoregressive process (BAR) when some of the data are missing. In this aim, we model the observed data by a two-type Galton-Watson process consistent with the binary Bee structure of the data. Under independence between the process leading to the missing data and the BAR process and suitable assumptions on the driven noise, we establish the strong consistency
of our estimators on the set of non-extinction of the Galton-Watson process, via a martingale approach. We also prove a quadratic strong law and the asymptotic normality.
\end{abstract}


\end{frontmatter}

\section{Introduction}
\label{section intro}
Bifurcating autoregressive processes (BAR)  generalize autoregressive (AR) processes, when the data have a binary tree structure. Typically, they are involved in modeling cell lineage data, since each cell in one generation gives birth to two offspring in the next one. Cell lineage data usually consist of observations of some quantitative characteristic of the cells, over several generations descended from an initial cell. BAR processes take into account both inherited and environmental effects to explain the evolution of the quantitative characteristic under study. They were first introduced by Cowan and Staudte \cite{CoSt86}. In their paper, the original BAR process was defined as follows. The initial cell is labelled $1$, and the two offspring of cell $k$ are labelled $2k$ and $2k+1$. If $X_k$ denotes the quantitative characteristic of individual $k$, then the first-order BAR process is given, for all $k\geq 1$, by
\begin{equation*}
\left\{\begin{array}{lcl}
X_{2k}&=&a+bX_k+\varepsilon_{2k},\\
X_{2k+1}&=&a+bX_k+\varepsilon_{2k+1}.\\
\end{array}\right.
\end{equation*}
The noise sequence $(\varepsilon_{2k}, \varepsilon_{2k+1})$ represents environmental effects, while $a,b$ are unknown real parameters, with $|b|<1$, related to the inherited effects. The driven noise $(\varepsilon_{2k}, \varepsilon_{2k+1})$ was originally supposed to be independent and identically distributed with normal distribution. But since two sister cells are in the same environment at their birth, $\varepsilon_{2k}$ and $\varepsilon_{2k+1}$ are allowed to be correlated, inducing a correlation between sister cells, distinct from the correlation inherited from their mother.

\smallskip

Recently, experiments made by biologists on aging of
\emph{Escherichia coli} \cite{SMT05}, motivated mathematical
and statistical studies of the asymmetric BAR process, that is when
the quantitative characteristics of the even and odd sisters are
allowed to depend on their mother's through different sets of
parameters $(a,b)$, see Equation~(\ref{defbar}) below. In
\cite{GBPSDT05, Guy07}, Guyon proposes an interpretation of the
asymmetric BAR process as a bifurcating Markov chain, which allows
him to derive laws of large numbers and central limit theorems for
the least squares estimators of the unknown parameters of the
process. This Markov chain approach was further developed by
Bansaye \cite{Ban08} in the context of cell division with parasite
infection, and by Delmas and Marsalle \cite{DM08},
where the cells are allowed to die. Another approach based on
martingales theory was proposed by Bercu, de Saporta and
G\'egout-Petit \cite{BSG09}, to sharpen the asymptotic analysis of
Guyon under weaker assumptions.

\smallskip

\begin{figure}[ht]
\centering
\includegraphics[height=5cm, angle=0]{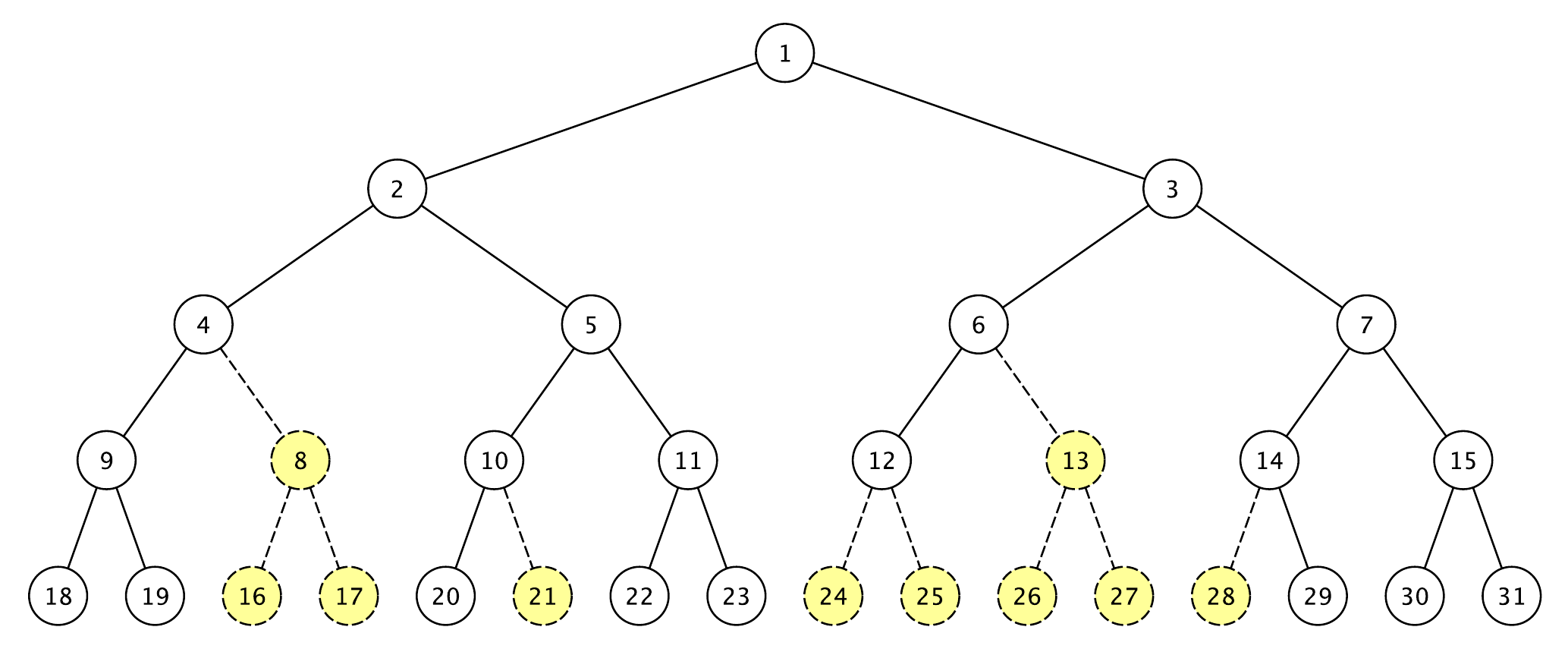}
\caption{A tree associated with the bifurcating auto-regressive process up to the 4th generation. The dashed cells are not observed.\label{arbre}}
\end{figure}

The originality of this paper is that we take into account possibly missing data in the estimation procedure of the parameters of the asymmetric BAR process, see Figure~\ref{arbre} for an example.
This is a problem of practical interest, as experimental data are often incomplete, either because some cells died, or because the measurement of the characteristic under study was impossible or faulty. For instance, among the 94 colonies dividing up to 9 times studied in \cite{SMT05}, in average, there are about 47\% of missing data. It is important to take this phenomenon into account in the model for a rigorous statistical study.

\smallskip

Missing data in bifurcating processes were first modeled by Delmas and Marsalle \cite{DM08}. They defined the genealogy of the cells through a Galton-Watson process, but they took into account the possible asymmetry problem only by differentiating the reproduction laws according to the daughter's type (even or odd). The bifurcating process was thus still a Markov chain. However, considering the biological issue of aging in \textit{E. coli} naturally leads to introduce the possibility that two cells of different types may not have the same reproduction law. In this paper, we thus introduce a two-type Galton-Watson process to model the genealogy, and lose the Markovian structure of the bifurcating chain, so that we cannot use the same approach as \cite{DM08}. Instead, we use the martingale approach introduced in \cite{BSG09}. It must be pointed out that missing data are not dealt with in \cite{BSG09}, so that we cannot directly use their results either. In particular, the observation process is another source of randomness that requires stronger moment assumptions on the driven noise of the BAR process and careful choice between various filtrations. In addition, the normalizing terms are now random and the convergences are only available on the random non-extinction set of the observed process.

\smallskip

The \emph{naive} approach to handle missing data would be to replace the sums over all data in the estimators by sums over the observed data only. Our approach is slightly more subtle, as we distinguish whether a cell has even or odd daughters. We propose a joint model where the structure for the observed data is based on a two-type Galton-Watson process consistent with the possibly asymmetric structure of the BAR process. See e.g. \cite{Har63, AN72, HJV07} for a presentation of multi-type Galton-Watson processes and general branching processes. Note also that our estimation procedure does not require the previous knowledge of the parameters of the two-type Galton-Watson process.

\smallskip

This paper is organized as follows. In Section ~\ref{section model}, we first introduce our BAR model as well as related notation, then we define and recall results on the two-type Galton-Watson process used to model the observation process. In Section~\ref{section LS}, we give the least square estimator for the parameters of observed BAR process and we state our main results on the convergence and asymptotic normality of our estimators as well as estimation results on data. The proofs are detailed in the following sections.

\section{Joint model}
\label{section model}
We now introduce our joint model, starting with the asymmetric BAR process for the variables of interest.
\subsection{Bifurcating autoregressive processes}
\begin{figure}[ht]
\centering
\includegraphics[height=7cm, angle=0]{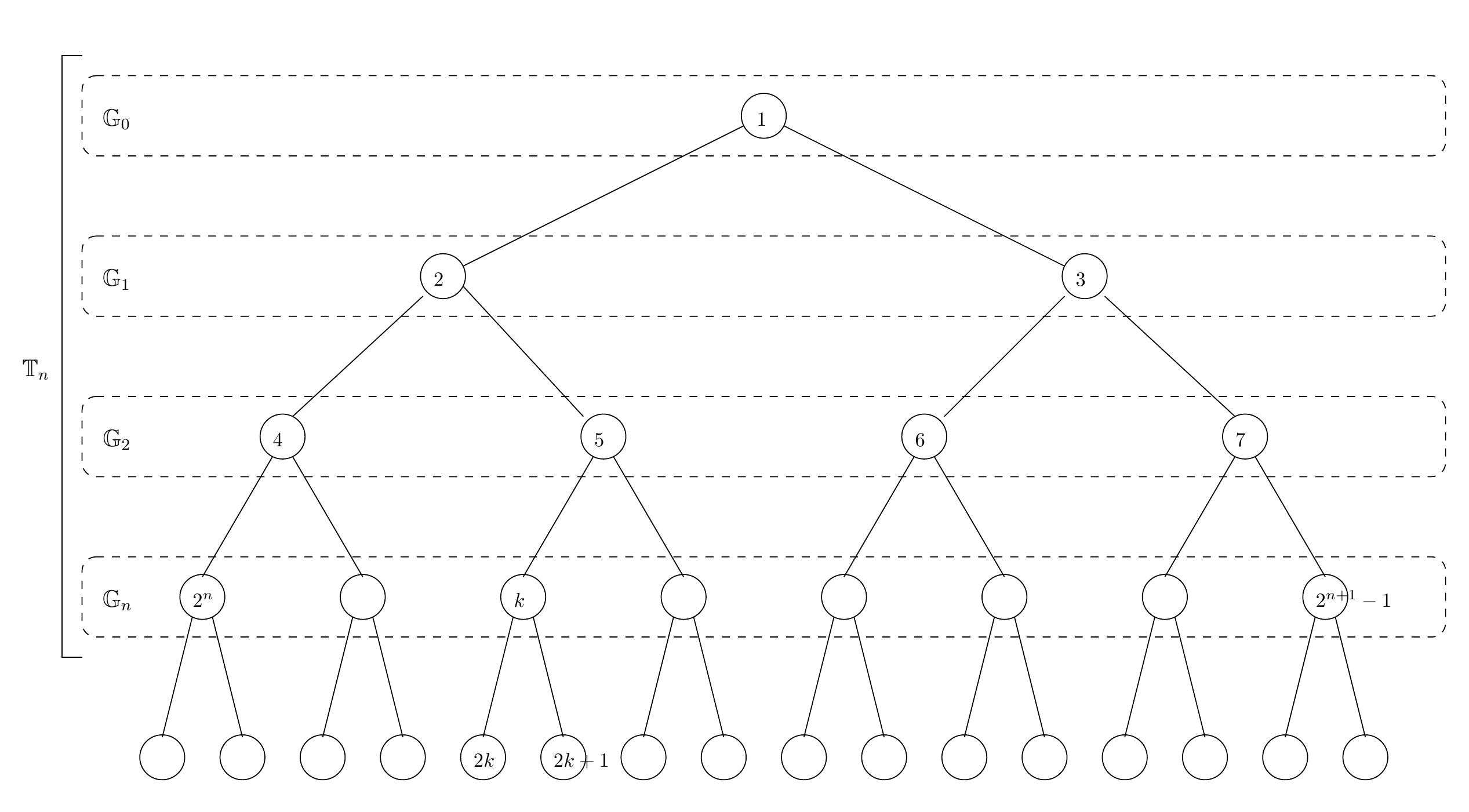}
\caption{The tree associated with the bifurcating auto-regressive process.}\label{FullTree}
\end{figure}

On the probability space $(\Omega, {\cal A }, \dP)$, we consider the first-order asymmetric BAR process given, for all $k\geq 1$, by
\begin{equation}\label{defbar}
\left\{
    \begin{array}{lcccccl}
     X_{2k} & = & a &+ &bX_k &+ &\varepsilon_{2k}, \\
     X_{2k+1} & = & c & +& dX_k &+ &\varepsilon_{2k+1}.
    \end{array}\right.
\end{equation}
The initial state $X_1$ is the characteristic of the ancestor, while $(\varepsilon_{2k},\varepsilon_{2k+1})$
is the driven noise of the process.
In all the sequel, we shall assume that $\mathbb{E}[X_1^8]<\infty$. Moreover, as in the previous literature, the parameters
$(a,b,c,d)$ belong to $\dR^4$ with
\begin{equation*}
  0<\max(|b|, |d|) < 1.
\end{equation*}
This assumption ensures the stability (non explosion) of the BAR process. As explained in the introduction, one can see this BAR process as a first-order autoregressive process on a binary tree,
where each vertex represents an individual or cell, vertex $1$ being the original ancestor, see Figure \ref{FullTree} for an illustration.
We use the same notation as in \cite{BSG09}. For all $n\geq 1$, denote the $n$-th generation by
$
 \dG_n = \{2^n, 2^n+1,\ldots,2^{n+1}-1\}$.
In particular, $\dG_0 = \{1\}$ is the initial generation, and $\dG_1 = \{2,3\}$ is the first generation of offspring from the first ancestor. Let $\dG_{r_k}$ be the generation of individual $k$, which means that $r_k=[\log_2(k)]$,  where $[x]$ denotes the largest integer less than or equal to $x$. Recall that the two offspring of individual $k$ are labelled $2k$ and $2k+1$, or conversely, the mother of individual $k$ is $[k/2]$.
More generally, the ancestors of individual $k$ are $[k/2], [k/2^2],\ldots,\ [k/2^{r_k}]$.
Denote by
 $\dT_n = \bigcup_{\ell=0}^n\dG_{\ell},
$
the sub-tree of all individuals from the original individual up to the $n$-th generation. Note that the cardinality $|\dG_n|$ of $\dG_n$  is $2^n$, while that of $\dT_n$ is $|\dT_n|=2^{n+1}-1$. Next, $\dT$ denotes the complete tree, so to speak $\dT=\bigcup_{n\ge0} \dG_n=\bigcup_{n\ge 0}\dT_n=\mathbb{N}^*=\dN\backslash\{0\}$. Finally, we need to distinguish the individuals in $\dG_n$ and $\dT_n$ according to their type. Since we are dealing with the types even and odd, that we will also label $0$ and $1$, we set 
\begin{equation*}
 \dG_n^0=\dG_n \cap (2 \mathbb{N}), \quad \dG_n^1=\dG_n \cap (2 \mathbb{N} + 1), \quad \dT_n^0=\dT_n \cap (2 \mathbb{N}),
\end{equation*}
\begin{equation} \label{def:arbretype}
\dT_n^1=\dT_n \cap (2 \mathbb{N}+1), \quad \dT^0=\dT \cap (2 \mathbb{N}) \quad \text{and} \quad \dT^1=\dT \cap (2 \mathbb{N}+1).
\end{equation}

\smallskip

We now state our assumptions on the noise sequence. Denote by $\dF=(\cF_n)$ the natural filtration associated with the
first-order BAR process, which means that $\cF_n$ is the $\sigma$-algebra generated by all individuals up to the
$n$-th generation, $\cF_n = \sigma\{X_k, k\in \dT_n\}$. In all the sequel, we shall make use of the following moment and independence hypotheses.
\begin{description}
\item[(HN.1)]  For all $n\geq 0$ and for all $k\in \dG_{n+1}$,
$\veps_k$ belongs to $L^8$. 
Moreover, there exist $(\sigma^2, \tau^4, \kappa^8) \in (0, +\infty)^3$, $(|\rho'|, \nu^2, \lambda^4) \in [0,1)^3$ such that :
\begin{itemize}
 \item $\forall n\geq 0 \text{ and }  k\in \dG_{n+1}$,
$$\dE[\veps_k|\cF_n] = 0, \hspace{0.1cm}  \dE[\veps_k^2|\cF_n]=\sigma^2, \hspace{0.1cm} \dE[\veps_k^4|\cF_n]=\tau^4, \hspace{0.1cm} \dE[\veps_k^8|\cF_n]=\kappa^8 \hspace{0.2cm}\text{a.s.}$$
 \item $\forall n\geq 0 \quad  \forall k\neq l \in \dG_{n+1}  \text{ with }[k/2]=[l/2]$,
$$\hspace{-1.8cm}\dE[\veps_k\veps_l|\cF_n] = \rho =\rho' \sigma^2, \hspace{0.1cm} \dE[\veps_{2k}^2\veps_{2k+1}^2|\cF_n] = \nu^2 \tau^4, \hspace{0.1cm} \dE[\veps_{2k}^4\veps_{2k+1}^4|\cF_n] = \lambda^4 \kappa^8   \hspace{0.2cm}\text{a.s.}$$
\end{itemize}

\item[(HN.2)] For all $n\geq 0$ the random vectors $\{(\veps_{2k} , \veps_{2k+1}), k \in\dG_n\}$ are conditionally independent given $\cF_n$.
\end{description}

\subsection{Observation process}
\label{section GW}
\begin{figure}[ht]
\centering
\includegraphics[height=6.9cm, angle=0]{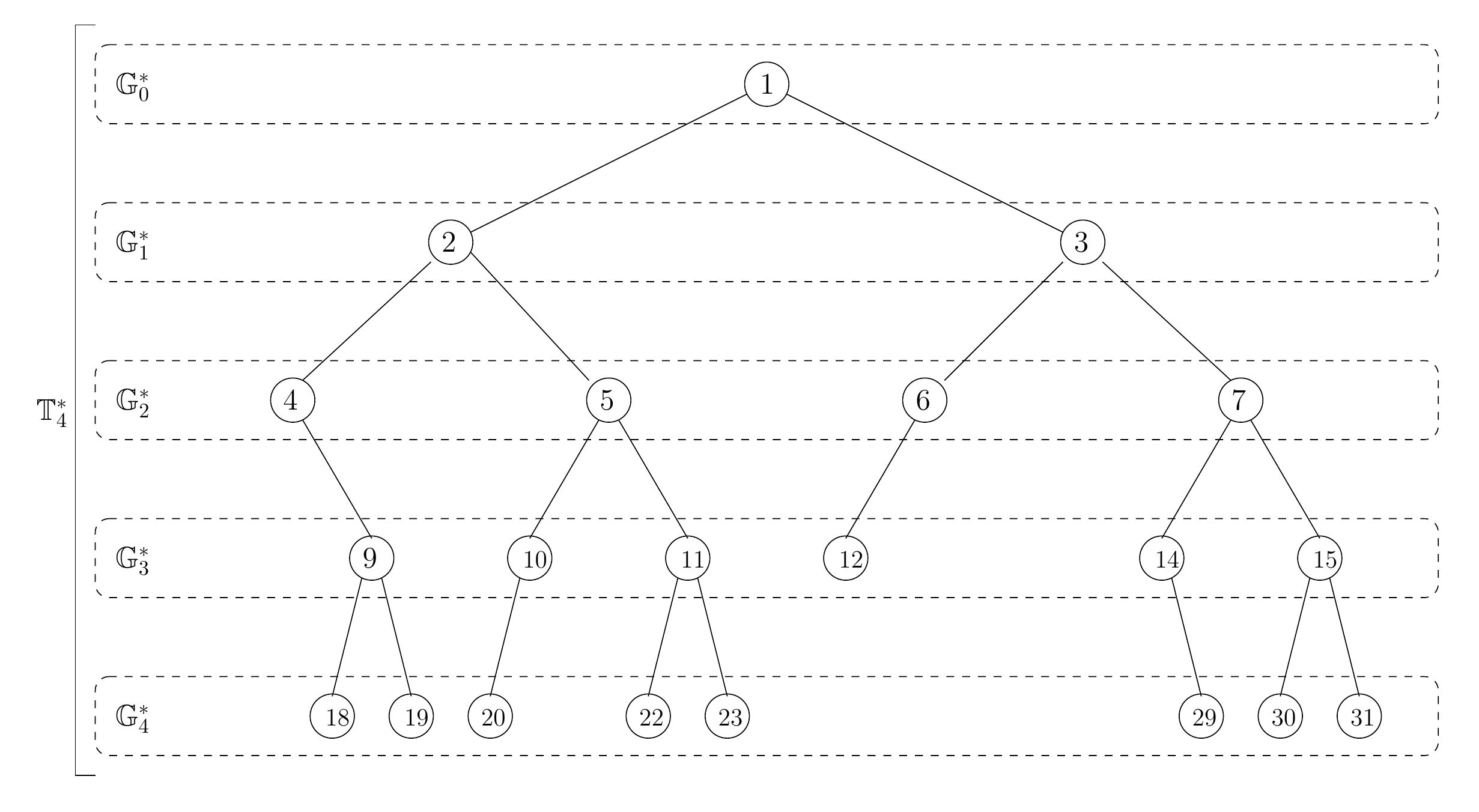}
\caption{The tree associated with the observed data of the tree in Figure~\ref{arbre}.}\label{MissingTree}
\end{figure}

We now turn to the modeling of the observation process. The observation process is intended to encode if a datum is missing or not. The natural property it has thus to satisfy is the following: if the datum is missing for some individual, it is also missing for all its descendants. Indeed, the datum may be missing because of the death of the individual, or because the individual is the last of its lineage at the end of the data's gathering, see Figure~\ref{MissingTree} for an example of partially observed tree.
\subsubsection{Definition of the observation process}
\label{section def GW}
Mathematically, we define the observation process, $(\delta_k)_{k\in\dT}$, as follows. We set $\delta_1=1$ and define recursively the sequence through the following equalities:
\begin{equation} \label{defdelta}
 \delta_{2k}=\delta_k \zeta_k^0 \quad \text{and} \quad \delta_{2k+1}=\delta_k \zeta_k^1,
\end{equation}
where $(\bs{\zeta}_k=(\zeta_k^0, \zeta_k^1))$ is a sequence of independent random vectors of $\{0,1\}^2$, $\zeta_k^i$ standing for the number ($0$ or $1$) of descendants of type $i$ of individual $k$. The sequences $(\bs{\zeta}_k, k \in 2\mathbb{N}^*)$ and $(\bs{\zeta}_k, k \in 2\mathbb{N}+1)$ are sequences of identically distributed random vectors. We specify the common laws of these two sequences using their generating functions, $f^{(0)}$ and $f^{(1)}$ respectively:
\begin{eqnarray*}
f^{(0)}(s_0, s_1)&=& p^{(0)}(0,0)+p^{(0)}(1,0)s_0+p^{(0)}(0,1)s_1 + p^{(0)}(1,1)s_0s_1,\\
f^{(1)}(s_0, s_1)&=& p^{(1)}(0,0)+p^{(1)}(1,0)s_0+p^{(1)}(0,1)s_1 + p^{(1)}(1,1)s_0s_1,
\end{eqnarray*}
where $p^{(i)}(j_0,j_1)$ is the probability that an individual of type $i$ gives birth to $j_0$ descendants of type $0$, and $j_1$ of type $1$. The sequence $(\delta_k)$ is thus completely defined. We also assume that the observation process is independent from the BAR process.
\begin{description}
 \item[(HI)] The sequences $(\delta_k)$ and $(\bs{\zeta}_k)$ are independent from the sequences $(X_k)$ and $(\veps_k)$.
\end{description}
Remark that, since both $\zeta^0_k$ and $\zeta^1_k$ take values in $\{0,1\}$ for all $k$, the observation process $(\delta_k)$ is itself taking values in $\{0,1\}$. Finally, Equation~(\ref{defdelta}) ensures that if $\delta_k=0$ for some $k \ge 2$, then for all its descendants $j$, $\delta_j=0$. In relation with the observation process $(\delta_k)$, we introduce two filtrations: $\cZ_n = \sigma \{ \bs{\zeta}_k , k\in \dT_n\}$, $\cO_n=\sigma\{\delta_k, k\in\dT_n\}$, and the sigma field $\cO = \sigma \{\delta_k , k\in \dT \}$. Notice that $\cO_{n+1} \subset \cZ_n$. We also define the sets of observed individuals as follows:
\begin{equation*}
\dG_n^*=\{k \in \dG_n : \delta_k =1\} \quad \text{and} \quad \dT_n^*=\{k \in \dT_n : \delta_k =1\}.
\end{equation*}
Finally, let $\mathcal{E}$ be the event corresponding to the cases when there are no individual left to observe. More precisely,
\begin{equation} \label{defext}
 \mathcal{E} = \bigcup_{n \ge 1} \{|\dG_n^*|=0 \}.
\end{equation}
We will denote $\overline{\cE}$ the complementary set of $\cE$.
\subsubsection{Results on the observation process}
\label{section res GW}
Let us introduce some additional notation. For $n \ge 1$, we define the number of observed individuals among the n-th generation, distinguishing according to their types:
\begin{equation}
 Z_n^0=|\dG_n^* \cap 2 \mathbb{N}| \quad \text{and} \quad Z_n^1=|\dG_n^* \cap (2 \mathbb{N}+1)|,
\end{equation}
and we set, for all $n \ge 1$, $\bs{Z}_n=(Z_n^0,Z_n^1)$. Note that for $i\in\{0,1\}$ and $n\geq 1$ one has
\begin{equation*}
 Z_n^i=\sum_{k\in\dG_{n-1}}\delta_{2k+i}.
\end{equation*}
One has $\dG_0^*=\dG_0=\{1\}$, but, even if $1$ is odd, the individual whose lineage we study may as well be of type $0$ as of type $1$. Consequently, we will work with possibly two different initial laws: $\dP(\cdot|\bs{Z}_0=\bs{e}_i)$, for $i \in \{0,1\}$, where $\bs{e}_0=(1,0)$ and $\bs{e}_1=(0,1)$. The process $(\bs{Z}_n, n \ge 0)$ is thus a two-type Galton-Watson process, and all the results we are giving in this section mainly come from \cite{Har63}. Notice that the law of $\bs{\zeta}_k$, for even $k$, is the law of reproduction of an individual of type $0$, the first component of $\bs{\zeta}_k$ giving the number of children of type $0$, the second the number of children of type $1$. The same holds for $\bs{\zeta}_k$, with odd $k$, \textit{mutatis mutandis}. This ensures the existence of moments of all order for these reproduction laws, and we can thus define the descendants matrix $\bs{P}$
\begin{equation*}
\bs{P}=\left(\begin{array}{cc}
           p_{00} & p_{01} \\
           p_{10} & p_{11}
          \end{array}\right),
\end{equation*}
where $p_{i0} = p^{(i)}(1,0) + p^{(i)}(1,1)$ and $p_{i1} = p^{(i)}(0,1) + p^{(i)}(1,1)$, for $i \in \{0,1\}$. The quantity $p_{ij}= \dE[\zeta_{2+i}^j]$ is thus the expected number of descendants of type $j$ of an individual of type $i$. We also introduce the variance of the laws of reproduction: $\sigma_{ij}^2 = \dE[(\zeta_{2+i}^j -p_{ij})^2]$, for $(i,j) \in \{0,1\}^2$. Note that $\sigma_{ij}^2 = p_{ij}(1-p_{ij})$. It is well-known (see e.g. Theorem 5.1 of \cite{Har63}) that when all the entries of the matrix $\bs{P}$ are positive, $\bs{P}$ has a positive strictly dominant eigenvalue, denoted $\pi$, which is also simple. We make the following main assumptions on the matrix $\bs{P}$.
\begin{description}
 \item[(HO)] All entries of the matrix $\bs{P}$ are positive: for all $(i,j) \in \{0,1\}^2$, $p_{ij} > 0$, and the dominant eigenvalue is greater than one: $\pi > 1$ .
\end{description}
Hence, still following Theorem 5.1 of \cite{Har63}, we know that there exist left and right eigenvectors for $\pi$ which are positive, in the sense that each component of the vector is positive. We call $\bs{y}=(y^0,y^1)^t$ such a right eigenvector, and $\bs{z}=(z^0,z^1)$ such a left one; without loss of generality, we choose $\bs{z}$ such that $z^0+z^1=1$. Regarding the two-type Galton-Watson process $(\bs{Z}_n)$, $\pi$ plays the same role as the expected number of offspring, in the case of standard Galton-Watson processes. In particular, $\pi$ is related to the extinction of the process, where the set of extinction of $(\bs{Z}_n)$ is defined as $\cup_{n \ge 1} \{\bs{Z}_n=(0,0)\}$. Notice that $\{\bs{Z}_n=(0,0)\} =\{Z_n^0 + Z_n^1=0\}= \{|\dG_n^*|=0 \}$, so that this set coincides with $\mathcal{E}$, defined by Eq.~\reff{defext}. Now let $\bs{q}=(q^0, q^1)$, where, for $i \in \{0,1\}$,
\begin{equation*}
q^i = \dP(\mathcal{E}|\bs{Z}_0=\bs{e}_i).
\end{equation*}
The probability $q^i$ is thus the extinction probability if initially there is one individual of type $i$. These two probabilities allow to compute the extinction probability under any initial distribution, since $\dP(\mathcal{E}) = \dE[(q^0)^{Z_0^0} (q^1)^{Z_0^1}]$, thanks to the branching property. Hypothesis (\textbf{HO}) means that the Galton-Watson process $(\bs{Z}_n)$ is super-critical, and ensures that $0 \le q^i <1$, for both $i=0$ and $i=1$. This immediately yields
\begin{equation}\label{probext}
\dP(\mathcal{E}) < 1.
\end{equation}
 Under that condition, we also have the existence of a non-negative random variable $W$ such that for any initial distribution of $\bs{Z}_0$
\begin{equation} \label{limGW}
\lim_{n \rightarrow +\infty} \frac{\bs{Z}_n}{\pi^n} = \lim_{n \rightarrow +\infty} \frac{\pi -1}{\pi^{n+1}-1}\sum_{\ell=0}^n \bs{Z}_{\ell} = W \bs{z} \quad \text{a.s.}.
\end{equation}
It is well-known that $\{W = 0\} = \mathcal{E}$ a.s., so that the set $\{W > 0\}$ can be viewed as the set of non-extinction $\overline{\cE}$ of $(\bs{Z}_n)$, up to a negligible set. These results give the asymptotic behavior of the number of observed individuals, since $|\dG^*_n|=Z_n^0 + Z_n^1$, and $|\dT^*_n|=\sum_{\ell=0}^n (Z_{\ell}^0 + Z_{\ell}^1$):
\begin{equation*}
 \lim_{n \rightarrow +\infty} \frac{|\dG^*_n|}{\pi^n} = \lim_{n \rightarrow +\infty} \frac{\pi -1}{\pi^{n+1}-1}|\dT^*_n| = W \hspace{0.5cm} \text{a.s.}
\end{equation*}
Roughly speaking, this means that $\pi^n$ is a deterministic equivalent of $|\dT^*_n|$ and Eq.~(\ref{limGW}) implies that $z^i$ is the asymptotic proportion of cells of type $i$ in a given generation. We will thus very often replace $|\dT^*_n|$ by $\pi^n$ for computations, and the next lemma will be used frequently to replace $\pi^n$ by $|\dT^*_n|$.
\begin{Lemma} \label{lemcard}
Under assumption \emph{(\textbf{HO})}, we have 
\begin{equation*}
 \lim_{n \rightarrow +\infty} \ind{|\dG_n^*|>0}  \frac{\pi^n}{|\dT^*_n|} = \frac{\pi -1}{\pi} \frac{1}{W } \indnex \quad \text{a.s.}
\end{equation*}
\end{Lemma}
\subsection{Joint model}
The model under study in this paper is therefore the observed BAR process defined by
\begin{equation*}
\left\{
    \begin{array}{rccccccl}
    \delta_{2k} X_{2k} & = & \delta_{2k}&(a &+ &bX_k &+ &\varepsilon_{2k}), \\
     \delta_{2k+1}X_{2k+1} & = &  \delta_{2k+1}&(c & +& dX_k &+ &\varepsilon_{2k+1}).
    \end{array}\right.
\end{equation*}
The aim of this paper is to study the sharp asymptotic properties of the least-squares estimators of the parameters $(a,b,c,d)$ and the variance matrix of the noise process.
\section{Least-squares estimation}
\label{section LS}
Our goal is to estimate ${\bs{\theta}}=(a,b,c,d)^t$ from the observed individuals up to the $n$-th generation, that is the
observed sub-tree $\dT_n^*$.
\subsection{Definition of the estimators}
\label{section def LS}
We propose to make use of the standard least-squares (LS) estimator $\wh{{\bs{\theta}}}_n$ which minimizes
\begin{equation*}
\Delta_n(\bs{\theta})=\sum_{k\in \dT_{n-1}}\delta_{2k}(X_{2k}-a-bX_k)^2 + \delta_{2k+1} (X_{2k+1}-c-dX_k)^2.
\end{equation*}
Consequently, we obviously have for all $n\geq 1$
\begin{equation}\label{defLS}
(\wh{\bs{\theta}}_n)= \left(
    \begin{array}{c}
     \wh{a}_n\\
     \wh{b}_n\\
     \wh{c}_n\\
     \wh{d}_n
    \end{array}\right) = \bs{\Sigma}_{n-1}^{-1}\sum_{k \in \mathbb{T}_{n-1}}\left(
\begin{array}{c}
\delta_{2k}X_{2k}  \\
\delta_{2k}X_kX_{2k} \\
\delta_{2k+1}X_{2k+1} \\
\delta_{2k+1}X_kX_{2k+1}
\end{array}\right),
\end{equation}
where, for all $n\geq 0$,
\begin{equation*}
\bs{\Sigma}_{n} = \left( \begin{array}{cc}
\bs{S}^0_{n} & 0 \\
0 & \bs{S}^1_{n}
\end{array} \right), \quad\textrm{and}\quad
\bs{S}^i_{n} =\sum_{k \in \mathbb{T}_{n}}\delta_{2k+i}\left(
\begin{array}{cc}
1 & X_k \\
X_k &X^2_k
\end{array}\right), 
\end{equation*}
for $i\in\{0,1\}$.
In order to avoid intricate invertibility assumption, we shall assume,
without loss of generality, that for all $n\geq 0$, $\bs{\Sigma}_n$ is invertible. Otherwise,
we only have to add the identity matrix $\rI_4$ to $\bs{\Sigma}_n$, as Proposition~\ref{mainlemma} states that the normalized limit of $\bs{\Sigma}_n$ is positive definite.
\begin{Remark}
Note that when all data are observed, that is when all $\delta_k$ equal $1$, this is simply the least squares estimator described in the previous literature. However, one must be careful here with the indices in the normalizing matrix, as there are now two different matrices $\bs{S}^0_{n}$ and $\bs{S}^1_{n}$, while there was only one in the fully observed problem. The intuitive way to deal with missing data would be to restrict the sums to the observed data only. Note that our estimator is more complex as it involves sums depending on the absence or presence of even- or odd-type \emph{daughters} of the available data.
\end{Remark}
We now turn to the estimation of the parameters $\sigma^2$ and $\rho$. We propose to estimate the conditional variance $\sigma^2$ and the conditional covariance $\rho$ by
\begin{equation*}
\wh{\sigma}^2_n = \frac{1}{ |\dT_{n}^*|}\sum_{k \in
\dT^*_{n-1}} (\wh{\veps}_{2k}^2 +   \wh{\varepsilon}_{2k+1}^2),\quad \wh{\rho}_n = \frac{1}{|\dT_{n-1}^{*01}|}\sum_{k \in\dT_{n-1}}
\wh{\veps}_{2k} \wh{\veps}_{2k+1},
\end{equation*}
where for all $k\in\mathbb{G}_n$,
\begin{equation*}
\left\{
    \begin{array}{lcrclcl}
     \wh{\veps}_{2k} &=& \delta_{2k}(X_{2k} &-& \wh{a}_n &-& \wh{b}_{n}X_k), \vspace{1ex}\\
     \wh{\veps}_{2k+1} &=& \delta_{2k+1}(X_{2k+1} &-& \wh{c}_{n} &-& \wh{d}_{n}X_k).
    \end{array}\right.,
\end{equation*}
and
\begin{equation*}
\dT_{n}^{*01} =\{k \in \dT_n : \delta_{2k}  \delta_{2k+1}=1\},
\end{equation*}
so to speak $\dT_{n-1}^{*01}$ is the set of the cells of the tree $\dT_{n-1}$ which have exactly two offspring.
\subsection{Main results}
\label{section results}
We can now state the sharp convergence results we obtain for the estimators above. We introduce additional notation.
For $i \in \{0,1\}$, let us denote :
\begin{equation*}
\bs{L}^i=
\left(\begin{array}{cc} \pi z^i & h^i \vspace{1ex}\\
h^i & k^i
\end{array}\right)\qquad \qquad \bs{L}^{0,1}=
\left(\begin{array}{cc}\bar{p}(1,1) & h^{0,1} \vspace{1ex}\\
h^{0,1} & k^{0,1}
\end{array}\right)
\end{equation*}
where $\bs{z}=(z^0,z^1)$ is the left eigenvector for the dominant eigenvalue $\pi$ of the descendants matrix $\bs{P}$ introduced in section \ref{section res GW}, $h^i$, $k^i$ are defined in Propositions~\ref{propSumX} and \ref{propSumX2} and the four terms of $L^{0,1}$ defined in Proposition~\ref{propS 01}. We also define the $4\times4$ matrices
\begin{equation}\label{def Lambda}
\bs{\Sigma} = \left(\begin{array}{cc}\bs{L}^0&0\\0&\bs{L}^1\end{array}\right), \quad\textrm{and}\quad\bs{\Gamma} = \left(\begin{array}{cc}\sigma^2\bs{L}^0&\rho\bs{L}^{0,1}\\\rho\bs{L}^{0,1}&\sigma^2\bs{L}^1\end{array}\right).
\end{equation}
Our first result deals with the strong consistency of the LS estimator $\wh{\bs{\theta}}_n$.
\begin{Theorem}\label{thmaptheta}
Under assumptions \emph{(\textbf{HN.1})}, \emph{(\textbf{HN.2})}, \emph{(\textbf{HO})} and \emph{(\textbf{HI})},
$\wh{\bs{\theta}}_{n}$ converges to $\bs{\theta}$ almost surely on $\overline{\cE}$ with the rate of convergence
\begin{equation}\label{thmaptheta1}
\ind{|\dG_n^*|>0}\| \widehat{\bs{\theta}}_{n}-\bs{\theta} \|^{2}=
\cO \left(\frac{\log |\dT_{n-1}^*|}{|\dT_{n-1}^*|} \right)\indnex
\hspace{1cm}\text{a.s.}
\end{equation}
In addition, we also have the quadratic strong law
\begin{equation}\label{thmaptheta2}
\lim_{n\rightarrow \infty}\ind{|\dG_n^*|>0} \frac{1}{n}\sum_{\ell=1}^n |\dT_{\ell-1}^*|
(\widehat{\bs{\theta}}_{\ell}-\bs{\theta})^t\bs{\Sigma} (\widehat{\bs{\theta}}_{\ell}-\bs{\theta})= 4\frac{\pi-1}{\pi}\sigma^2\indnex \hspace{1cm} \text{a.s.}
\end{equation}
\end{Theorem}
Our second result is devoted to the almost sure asymptotic properties of the variance and covariance estimators $\wh{\sigma}^2_n$ and $\wh{\rho}_n$. Let
\begin{equation*}
\sigma^2_n =\frac{1}{ |\dT_{n}^*|}\sum_{k
\in\dT^*_{n-1}}\!\!\!\!(\delta_{2k}{\veps}_{2k}^2 +
\delta_{2k+1}{\varepsilon}_{2k+1}^2),\quad
{\rho}_n = \frac{1}{|\dT_{n-1}^{*01}|}\sum_{k \in\dT^*_{n-1}}\!\!\!\!{\delta_{2k}\veps}_{2k}\delta_{2k+1}{\veps}_{2k+1}.
\end{equation*}
\begin{Theorem}\label{thmapsigmarho}
Under assumptions \emph{(\textbf{HN.1})}, \emph{(\textbf{HN.2})}, \emph{(\textbf{HO})} and \emph{(\textbf{HI})}, $\wh{\sigma}^2_n$ converges almost surely to $\sigma^2$ on $\overline{\cE}$. More precisely, one has
\begin{equation}\label{apsigma1}
\lim_{n\rightarrow\infty}\ind{|\dG_n^*|>0}\frac{1}{n}\sum_{k\in\dT_{n-1}}\sum_{i=0}^{1}\delta_{2k+i}(\wh{\varepsilon}_{2k+i}-\varepsilon_{2k+i})^2={4}{(\pi-1)}\sigma^2\indnex\hspace{.2cm}\text{a.s.}
\end{equation}
\begin{equation}\label{apsigma2}
\lim_{n\rightarrow\infty}\ind{|\dG_n^*|>0}\frac{|\dT_{n}^*|}{n}
(\wh{\sigma}^2_n-\sigma_n^2)
={4}{(\pi-1)}\sigma^2\indnex\hspace{1cm}\text{a.s.}
\end{equation}
In addition, $\wh{\rho}_n$ converges almost surely to $\rho$ on $\overline{\cE}$ and one has
\begin{align}
\lim_{n\rightarrow\infty}\ind{|\dG_n^*|>0}&\frac{1}{n}\sum_{k\in\dT_{n-1}}\delta_{2k}(\wh{\varepsilon}_{2k}-\varepsilon_{2k})\delta_{2k+1}(\wh{\varepsilon}_{2k+1}-\varepsilon_{2k+1})\nonumber\\
&=4\rho\frac{\pi-1}{\pi}tr\big((\bs{L}^1)^{-1/2}\bs{L}^{0,1}(\bs{L}^0)^{-1/2}\big)\indnex\hspace{0.5cm}\text{a.s.}\label{aprho1}
\end{align}
\begin{equation}\label{aprho2}
\lim_{n\rightarrow\infty}\ind{|\dG_n^*|>0}\frac{|\dT_{n}^*|}{n}
(\wh{\rho}_n-\rho_n)=\rho\frac{\pi-1}{\bar{p}(1,1)}tr\big((\bs{L}^1)^{-1}(\bs{L}^{0,1})^2(\bs{L}^0)^{-1}\big)\indnex
\hspace{0.5cm}\text{a.s.}
\end{equation}
\end{Theorem}
\noindent
Our third result concerns the asymptotic normality for all our estimators
$\wh{\bs{\theta}}_n$, $\wh{\sigma}^2_n$ and $\wh{\rho}_n$ given the non-extinction of the underlying Galton-Watson process.  For this, using the fact that $\dP(\overline{\cE})\neq 0$ thanks to Eq. \reff{probext}, we define the probability $ \dP_{\overline{\cE}}$ by
$$ \dP_{\overline{\cE}}(A) =  \frac{\dP(A \cap \overline{\cE})}{ \dP(\overline{\cE})}\qquad \text{ for all } A \in {\cal A}.$$
\begin{Theorem}\label{thmCLT}
Under assumptions \emph{(\textbf{HN.1})}, \emph{(\textbf{HN.2})}, \emph{(\textbf{HO})} and \emph{(\textbf{HI})}, we have the central limit theorem
\begin{equation}\label{CLTtheta}
\sqrt{|\dT^*_{n-1}|} (\widehat{\bs{\theta}}_{n}-\bs{\theta})
\liml
\cN(0,\bs{\Sigma}^{-1}\bs{\Gamma} \bs{\Sigma}^{-1})\quad\text{on }\ 
({\overline{\cE}}, \dP_{\overline{\cE}}).
\end{equation}
In addition, we also have
\begin{equation}
\label{CLTsigma}
\sqrt{|\dT^*_{n}|} (\wh{\sigma}^2_n-{\sigma}^2)
\liml
\cN\Bigl(0,\frac{\pi(\tau^4-\sigma^4)+2\bar{p}(1,1)(\nu^2 \tau^4-\sigma^4)}{\pi}\Bigr)\quad\text{on }\ 
({\overline{\cE}}, \dP_{\overline{\cE}}),
\end{equation}
where $\bar{p}(1,1)$ is defined in Eq. \reff{def pbar} and
\begin{equation}
\label{CLTrho}
\sqrt{|\dT^{*01}_{n-1}|} (\wh{\rho}_n-{\rho})
\liml
\cN(0,{\nu^2\tau^4-\rho^2})\quad\text{on }\ 
({\overline{\cE}}, \dP_{\overline{\cE}}).
\end{equation}
\end{Theorem}

The proof of our main results is going to be detailed in the next sections. It is based on martingale properties, and we will exhibit our main martingale $(\bs{M}_n)$ in Section~\ref{sectionmartingale}. Sections~\ref{LLN} to \ref{section proof M} are devoted proving to the sharp asymptotic properties of $(\bs{M}_n)$. Finally, in Section~\ref{section main proof} we prove our main results. Before turning to the definition of the martingale $(\bs{M}_n)$, we present a short application of our estimation procedure on data.

\subsection{Results  on real data}
\label{sectiondata}
The biological issue  addressed by Stewart et al. in \cite{SMT05} is aging in the single cell organism \emph{Escherichia coli}, see also \cite{ES07} for further biological details. \emph{E. coli} is a rod-shaped bacterium that reproduces by dividing in the middle. Each cell has thus a new end (or \emph{pole}), and an older one. The cell that inherits the old pole of its mother is called the old pole cell, the cell that inherits the new pole of its mother is called the new pole cell. Therefore, each cell has a \emph{type}: old pole (even) or new pole (odd) cell, inducing asymmetry in the cell division.

\medskip

Stewart et al. filmed colonies of dividing cells, determining the complete lineage and the growth rate of each cell. Their statistical study of the averaged genealogy and pair-wise comparison of sister cells showed that the old pole cells exhibit cumulatively slowed growth, less offspring biomass production and an increased probability of death. Note that their test assumes independence between the averaged pairs of sister cells which is not verified in the lineage. 

\medskip

Another analysis was proposed in \cite{GBPSDT05}. They model the growth rate by a Markovian bifurcating process, allowing single-experiment statistical analysis instead of averaging all the genealogical trees. Asymptotic properties of a more general asymmetric Markovian bifurcating autoregressive process are then investigated in \cite{Guy07}, where a Wald's type test is rigorously constructed to study the asymmetry of the process. These results cannot be compared to ours because this model does not take into account the possibly missing data from the genealogies, and it is not clear how the author manages them, as not a single tree from the data of \cite{SMT05} is complete.  In \cite{DM08}, the authors take missing data into account but, contrary to our approach, they  allow different sets of parameters for cells with two, one or no offspring, making the direct comparison with our estimator again impossible.
 
\medskip 

We have applied our methodology on the set of data \verb+penna-2002-10-04-4+ from the experiments of \cite{SMT05}. 
It is the largest data set of the experiment. It contains 663 cells up to generation 9 (note that there would be 1023 cells in a full tree up to generation 9).  In particular, we have performed
\begin{itemize}
\item point estimation of the vector $\bs{\theta}$,
\item interval estimation for the coefficients $(a,b,c,d)$,
\item Wald's type symmetry tests for the entries of $\wh{\bs{\theta}}_n$.
\end{itemize} 
\begin{table}[htdp]
\begin{center}
\begin{tabular}{|c|c|c|}
\hline
parameter & $a$ & $c$ \\
\hline
estimation   & 0.03627 &  0.03058   \\
C.I. & $[0.03276 ; 0.03979]$ &$[ 0.02696 ; 0.03420]$ \\
\hline
parameter & b & d \\
\hline
 estimation & 0.02662    &  0.17055 \\
C.I. &  $[- 0.06866  ; 0.12191]$ & $[0.07247 ; 0.26863]$  \\
 \hline
\end{tabular}
\caption{Estimation on the data set \emph{penna-2002-10-04-4} }\label{default}
\end{center}
\end{table}
Table 1 gives the estimation $\wh{\bs{\theta}}_9$ of $\wh{\bs{\theta}}$ with the 95\%  Confidence Interval (C.I.) of each coefficient. The variance given by the CLT for $\bs{\theta}$ in Eq.~\reff{CLTtheta}, is approximated by $\bs{\Sigma}^{-1}_n\bs{\Gamma}_n \bs{\Sigma}^{-1}_n$ thanks to the convergence given in Corollary~\ref{cv Gamman}. The confidence intervals of $b$ and $d$ show that the \emph{non explosion} assumption ($|b|<1$ and $|d|<1$) is satisfied. Some empirical computation on the process $(\delta_k)$ gives the following estimation for the highest eigenvalue  of the Galton-Watson process : $\hat{\pi}= 1.35669$ (with confidence interval $[1.27979, 1.43361]$, see \cite{MaTou00}), also satisfying the super-criticality assumption. Wald tests of comparison between the coefficients of $\bs{\theta}$ have been deduced of the CLT. The null hypotheses $(a,b)=(c,d)$ (resp. $a =c$, $b=d$) are rejected with p-values p= 0.0211 (resp. p= 0.0158 and p=0.0244). Hence on this data set the cell division is indeed statistically asymmetric.
\section{Martingale approach}
\label{sectionmartingale}
To establish all the asymptotic properties of our estimators, we shall make use of a martingale approach, similar to \cite{BSG09}. However, their results cannot be used in our framework, since the randomness comes now not only from the state process, but also from the time space (genealogy). These two mixed randomness sources require careful choice between various filtrations, and stronger moment assumptions on the driven noise of the BAR process. For all $n\geq 1$, denote
\begin{equation*}
\bs{M}_n= \sum_{k \in \dT_{n-1}}
\left(\delta_{2k}\veps_{2k},\  
\delta_{2k}X_k\veps_{2k},\
\delta_{2k+1}\veps_{2k+1},\
\delta_{2k+1}X_k\veps_{2k+1}\right)^t.
\end{equation*}
Thus, for all $n\geq2$, we readily deduce from Equations~(\ref{defLS}) and (\ref{defbar}) that
\begin{equation}\label{thetadiff}
\wh{\bs{\theta}}_n-\bs{\theta} =\bs{\Sigma}^{-1}_{n-1}
\sum_{k \in \dT_{n-1}}
\left( \begin{array}{cccc}
\delta_{2k}\veps_{2k}  \\
\delta_{2k}X_k\veps_{2k} \\
\delta_{2k+1}\veps_{2k+1} \\
\delta_{2k+1}X_k\veps_{2k+1}
\end{array}\right)=\bs{\Sigma}^{-1}_{n-1}\bs{M}_{n}.
\end{equation}
The key point of our approach is that $(\bs{M}_n)$ is a martingale for a well chosen filtration.
\subsection{Martingale property}
\label{section prop mart}
Recall that $\cO=\sigma\{\delta_k, k\in\dT\}$ is the $\sigma$-field generated by the observation process. We shall assume that all the history of the process $(\delta_k)$ is known at time $0$ and use the filtration $\dF^{\cO}=(\cF^{\cO}_n)$ defined for all $n$ by
\begin{equation*}
\cF^{\cO}_n=\cO\vee\sigma\{\delta_kX_k, k\in\dT_n\}=\cO\vee\sigma\{X_k, k\in\dT_n^*\},
\end{equation*}
where $\cF\vee\cG$ denotes the $\sigma$-field generated by both $\cF$ and $\cG$.
Note that for all $n$, $\cF^{\cO}_n$ is a sub $\sigma$-field of $\cO\vee\cF_n$.
\begin{Proposition}\label{prop M mart}
Under assumptions \emph{(\textbf{HN.1})}, \emph{(\textbf{HN.2})} and  \emph{(\textbf{HI})}, the process $(\bs{M}_n)$ is a square integrable $\dF^{\cO}$-martingale with increasing process given, for all $n\geq 1$, by
\begin{equation*}
<\bs{M}>_n = \bs{\Gamma}_{n-1} =  \left( \begin{array}{cc}
 \sigma^2\bs{S}^0_{n-1} & \rho \bs{S}^{0,1}_{n-1} \\
\rho \bs{S}^{0,1}_{n-1} & \sigma^2 \bs{S}^1_{n-1}
\end{array}\right),
\end{equation*}
where $\bs{S}^{0}_{n}$ and $\bs{S}^{1}_{n}$ are defined in section~\ref{section def LS} and
\begin{equation*}
\bs{S}^{0,1}_{n} =\sum_{k \in \mathbb{T}_{n}}\delta_{2k}\delta_{2k+1}
\left( \begin{array}{cc}
1 & X_k \\
X_k & X^2_k
\end{array}\right).
\end{equation*}
\end{Proposition}

\noindent\textbf{Proof :} First, notice that for all $n\geq1$, one has
\begin{equation*}
\Delta \bs{M}_n=\bs{M}_n-\bs{M}_{n-1}=\sum_{k \in \dG_{n-1}}
\left( \begin{array}{cccc}
\delta_{2k}\veps_{2k}  \\
\delta_{2k}X_k\veps_{2k} \\
\delta_{2k+1}\veps_{2k+1} \\
\delta_{2k+1}X_k\veps_{2k+1}
\end{array}\right).
\end{equation*}
Now, we use the fact that for all $n$, $\cF^{\cO}_n$ is a sub-$\sigma$ field of $\cO\vee\cF_n$, the independence between $\cO$ and $\cF_n$ under assumption (\textbf{HI}) and the moment hypothesis (\textbf{HN.1}) to obtain
\begin{eqnarray}\label{esp cond}
\dE[\delta_{2k}\veps_{2k}\ |\ \cF^{\cO}_{n-1}]
&=&\delta_{2k}\dE\big[\dE[\veps_{2k}\ |\ \cO\vee\cF_{n-1}]\ |\ \cF^{\cO}_{n-1}\big]\nonumber\\
&=&\delta_{2k}\dE\big[\dE[\veps_{2k}\ |\ \cF_{n-1}]\ |\ \cF^{\cO}_{n-1}\big]\ =\ 0.\nonumber
\end{eqnarray}
We obtain similar results for the other entries of $\Delta \bs{M}_n$ as $\delta_{2k+1}$ and $X_k$ are also $\cF_{n-1}^{\cO}$-measurable. Hence, $(\bs{M}_n)$ is a $\dF^{\cO}$-martingale. It is clearly square integrable from assumption (\textbf{HN.1}). The same measurability arguments together with assumption (\textbf{HN.2}) yield
\begin{eqnarray*}
\lefteqn{\dE[\Delta \bs{M}_n(\Delta \bs{M}_n)^t\ |\ \cF^{\cO}_{n-1}]}\\
&=&\!\!\!\!\sum_{k\in \dG_{n-1}}\left(\begin{array}{cccc}
\sigma^2\delta_{2k} &   \sigma^2\delta_{2k}X_k & \rho  \delta_{2k}\delta_{2k+1} &  \rho  \delta_{2k}\delta_{2k+1}X_k  \\
 \sigma^2\delta_{2k}X_k &   \sigma^2\delta_{2k}X_k^2 & \rho  \delta_{2k}\delta_{2k+1}X_k &  \rho  \delta_{2k}\delta_{2k+1}X_k^2\\
\rho  \delta_{2k}\delta_{2k+1} &  \rho  \delta_{2k}\delta_{2k+1}X_k &  \sigma^2\delta_{2k+1} &   \sigma^2\delta_{2k+1}X_k\\
\rho  \delta_{2k}\delta_{2k+1}X_k &  \rho  \delta_{2k}\delta_{2k+1}X_k^2 &   \sigma^2\delta_{2k+1}X_k &   \sigma^2\delta_{2k+1}X_k^2
\end{array}\right).
\end{eqnarray*}
Hence the result as $<\bs{M}>_n=\sum_{\ell=1}^{n}\dE[\Delta \bs{M}_{\ell}(\Delta \bs{M}_{\ell})^t\ |\ \cF^{\cO}_{\ell-1}]$.\hspace{\stretch{1}}$ \Box$\\

Our main results are direct consequences of the sharp asymptotic properties of the martingale $(\bs{M}_n)$. In particular, we will extensively use the strong law of large numbers for locally square integrable real martingales given in Theorem 1.3.15 of \cite{Duflo97}. Throughout this paper, we shall also use other auxiliary martingales, either with respect to the same filtration $\dF^{\cO}$, or with respect to other filtrations naturally embedded in our process, see Lemma~\ref{lemcvmart}.
\subsection{Asymptotic results}
\label{section conv mart}
We first give the asymptotic behavior of the matrices $\bs{S}^{0}_n$, $\bs{S}^1_n$ and $\bs{S}^{0,1}_n$. This is the first step of our asymptotic results.
\begin{Proposition}
\label{mainlemma}
Suppose that assumptions~\emph{(\textbf{HN.1})}, \emph{(\textbf{HN.2})}, \emph{(\textbf{HO})} and \emph{(\textbf{HI})} are satisfied. Then, for $i\in\{0,1\}$, we have
\begin{equation*}
\lim_{n\rightarrow \infty} \ind{|\dG_n^*|>0}\frac{\bs{S}^i_{n}}{|\dT_n^*|} = \indnex\bs{L}^i\hspace{0.3cm} \text{a.s.}\quad\text{and}\quad
 \lim_{n\rightarrow \infty} \ind{|\dG_n^*|>0}\frac{\bs{S}^{0,1}_{n}}{|\dT_n^*|} = \indnex\bs{L}^{0,1} \hspace{0.3cm} \text{a.s.}
\end{equation*}
In addition, $\bs{L}^0$ and $\bs{L}^1$, hence $\bs{\Sigma}$ are definite positive.
\end{Proposition}
A consequence of this proposition is the asymptotic behavior of the increasing process of the martingale $(\bs{M}_n)$.
\begin{Corollary}\label{cv Gamman}
Suppose that assumptions~\emph{(\textbf{HN.1})}, \emph{(\textbf{HN.2})}, \emph{(\textbf{HO})} and \emph{(\textbf{HI})} are satisfied. Then, we have
\begin{equation*}
\lim_{n\rightarrow \infty} \ind{|\dG_n^*|>0} \frac{\bs{\Sigma}_n}{|\dT_n^*|} =\indnex\bs{\Sigma},\quad
\text{and}\quad
\lim_{n\rightarrow \infty} \ind{|\dG_n^*|>0} \frac{\bs{\Gamma}_n}{|\dT_n^*|} =\indnex\bs{\Gamma}.
\end{equation*}
\end{Corollary}
This result is the keystone of our asymptotic analysis. It enables us to prove sharp asymptotic properties for the martingale $(\bs{M}_n)$.
\begin{Theorem}\label{th Mn}
Under assumptions \emph{(\textbf{HN.1})}, \emph{(\textbf{HN.2})}, \emph{(\textbf{HO})} and \emph{(\textbf{HI})}, we have
\begin{equation}\label{th Mn1}
\ind{|\dG_n^*|>0}\bs{M}_n^t\bs{\Sigma}_{n-1}^{-1}\bs{M}_n=\cO(n)\hspace{1cm} \text{a.s.}
\end{equation}
In addition, we also have
\begin{equation}\label{th Mn3}
\lim_{n\rightarrow\infty}\ind{|\dG_n^*|>0}\frac{1}{n}\sum_{\ell=1}^n \bs{M}_{\ell}^{t}\bs{\Sigma}_{\ell-1}^{-1}\bs{M}_{\ell}=4\sigma^2\indnex\hspace{1cm} \text{a.s.}
\end{equation}
Moreover, 
we have the central limit theorem on  $({\overline{\cE}}, \dP_{\overline{\cE}})$
\begin{equation*}
\frac{1}{\sqrt{|\dT^*_{n-1}|}}\bs{M}_n
\liml
\cN(0,\bs{\Gamma})\quad\text{on }\ 
({\overline{\cE}}, \dP_{\overline{\cE}}).
\end{equation*}
\end{Theorem}
As seen in Eq.~\reff{thetadiff}, $(\wh{\bs{\theta}}_n-\bs{\theta})$ is closely linked to $\bs{M}_n$ and this last theorem is then the major step to establish the asymptotic properties of our estimators.  The proof of this Theorem is given in Section~\ref{section proof M}. As explained before, it is a consequence of Proposition~\ref{mainlemma} which proof is detailed in Section~\ref{section proof mainlemma}. In between, Section~\ref{LLN} presents preliminary results in the form of laws of large number for the observation, noise and BAR processes.

\section{Laws of large numbers}
\label{LLN}
We now state some laws of large numbers involving the observation, noise and BAR processes. They are based on martingale convergence results, and we start with giving a general result of convergence for martingales adapted to our framework.
\subsection{Martingale convergence results}
\label{section auxiliary}
The following result is nothing but the strong law of large numbers for square integrable martingales, written in our peculiar setting, and will be repeatedly used.
\begin{Lemma}\label{lemcvmart}
Let $\mathcal{G}=(\cG_n)$ be some filtration, $(H_n)$ and $(G_n)$ be two sequences of random variables satisfying the following hypotheses:
\begin{description}
 \item[(i)] for all $n \ge 1$, for all $k \in \dG_n$, $H_k$ is $\cG_{n-1}$-measurable, $G_k$ is $\cG_n$-measurable, and $\dE[(H_kG_k)^2]< +\infty$,
 \item[(ii)] there exist $c^2>0$, $r \in [-1,1]$, such that for all $n \ge 1$, for all $k,p \in \dG_n$,
$$\dE[G_k | \cG_{n-1}]=0, \quad \dE[G_kG_p | \cG_{n-1}] = \left\{ \begin{array}{ll}
                                     c^2 & \text{if $k=p$}, \\
                                     rc^2 & \text{if $k \neq p$ and $[k/2]= [p/2]$}, \\
                                     0 & \text{otherwise,}
                                    \end{array} \right.$$
 \item[(iii)] there exists a sequence of real numbers $(a_n)$ that tends to $\infty$ such that $\sum_{k \in \dT_n} H_k^2 = \cO(a_n).$ 
\end{description}
Then $\sum_{k \in \dT_n} H_k G_k$ is a $\cG$-martingale and one has
$$\lim_{n \rightarrow \infty} \frac{1}{a_n} \sum_{k \in \dT_n} H_k G_k = 0 \hspace{0.5cm} a.s.$$
\end{Lemma}

\noindent\textbf{Proof:} Define $D_n = \sum_{k \in \dT_n} H_k G_k$. Assumptions (i) and (ii) clearly yield that $(D_n)$ is a square integrable martingale with respect to the filtration $(\cG_n)$. Thanks to (ii), its increasing process satisfies
\begin{eqnarray*}
 <D>_n &=& c^2  \big(\sum_{k \in \dT_n} H_k^2 +  2r\sum_{k \in \dT_{n-1}} H_{2k}H_{2k+1}\big) \\
       &\le& c^2  \big(\sum_{k \in \dT_n} H_k^2 +  r\sum_{k \in \dT_{n-1}} (H_{2k}^2+H_{2k+1}^2)\big)\\
       &\le& c^2 (r+1) \sum_{k \in \dT_n} H_k^2,
\end{eqnarray*}
and now, (iii) implies that $<D>_n= \cO(a_n)$. Finally, since the sequence $(a_n)$ tends to $\infty$, Theorem 1.3.15 of \cite{Duflo97} ensures that $D_n=o(a_n)$ a.s.
\hspace{\stretch{1}}$\Box$\\

We also recall Lemma~A.3 of \cite{BSG09} that will be useful in the sequel.
\begin{Lemma}\label{lemtoepbar}
Let $(\bs{A}_n)$ be a sequence of real-valued matrices such that
\begin{equation*}
 \sum_{n=0}^{\infty}\|\bs{A}_n\|<\infty\qquad\textrm{and}\qquad\lim_{n \rightarrow \infty} \sum_{k=0}^n\bs{A}_k=\bs{A}.
\end{equation*}
In addition, let $(\bs{X}_n)$ be a sequence of real-valued vectors which converges to a limiting value $\bs{X}$.
Then,
\begin{equation*}
\lim_{n \rightarrow \infty}{\sum_{\ell=0}^n \bs{A}_{n-\ell}\bs{X}_{\ell}}=\bs{AX}.
\end{equation*}
\end{Lemma}

\subsection{Laws of large numbers for the observation process}
\label{section LLN obs}
We now give more specific results on the asymptotic behavior of the observation process $(\delta_k)_{k \geq 1}$. Recall the notation $\dT_n^i$ defined in (\ref{def:arbretype}).

\begin{Lemma}\label{lemSumdelta}
Under the assumption \emph{(\textbf{HO})}, we have the
following convergences, for $(i,j)$ in $\{0,1\}^2$
\begin{equation*}
\lim_{n\rightarrow +\infty} \frac{1}{\pi^n}\sum_{k \in \dT_n^i}\delta_{2k+j} = p_{ij} \frac{\pi}{\pi-1}Wz^i \hspace{1cm}\text{a.s.}
\end{equation*}
\begin{equation*}
\lim_{n\rightarrow +\infty} \frac{1}{\pi^n}\sum_{k \in \dT_n^i}\delta_{2k}\delta_{2k+1} = p^{(i)}(1,1) \frac{\pi}{\pi-1}Wz^i \hspace{1cm}\text{a.s.}
\end{equation*}
\end{Lemma}

\noindent\textbf{Proof:} Recall that $\delta_{2k+j}=\delta_{k} \zeta^j_k$, so that
\begin{eqnarray*}
\sum_{k \in \dT_n^i}\delta_{2k+j}
& = & p_{ij} \sum_{k \in \dT_n^i}\delta_k +  \sum_{k \in \dT_n^i}\delta_k(\zeta^j_k - p_{ij})\ = \ p_{ij} \Big(i+\sum_{\ell=1}^n Z_{\ell}^i\Big) + D_n,
\end{eqnarray*}
since $\dG_0=\{1\}$, so that $\dT_n^i$ contains $1$ or not, according to $i=1$ or not, and where $D_n =\sum_{k \in \dT_n^i}\delta_k(\zeta^j_k - p_{ij})$. To deal with $D_n$, we use Lemma~\ref{lemcvmart}, with $\cG=(\cZ_n)$ (recall that $\cZ_n = \sigma \{ \bs{\zeta}_k , k\in \dT_n\}$), $H_k= \delta_k \ind{k \in \dT^i}$, and $G_k =(\zeta^j_k - p_{ij})\ind{k \in \dT^i}$. Assumption (i) of Lemma~\ref{lemcvmart} is obviously satisfied, since $\delta_k$, for $k \in \dG_n$, is $\cZ_{n-1}$-measurable. Regarding (ii), since the sequence $(\zeta^j_k)$ is a sequence of i.i.d. random variables with expectation $p_{ij}$ and variance $\sigma_{ij}^2$, we have $\dE[G_k | \cZ_{n-1}]=0$ and $\dE[G_k^2 | \cZ_{n-1}]=\sigma_{ij}^2$, for $k \in \dG_n$, and  $\dE[G_k G_p| \cZ_{n-1}]=0$, for $k \neq p \in \dG_n$. Finally, we turn to assumption (iii):
\begin{equation*}
\sum_{k \in \dT_n} H_k^2 = \sum_{k \in \dT_n^{i}}\delta_k = i+\sum_{\ell=1}^n Z_{\ell}^i = \cO(\pi^n),
\end{equation*}
thanks to \textbf{(HO)} and Eq.~\reff{limGW}. Finally, $D_n = o(\pi^n)$, and again using Eq.~\reff{limGW}, we obtain the first limit.
The proof of the second one is similar using the $\cZ$-martingale:
\begin{equation*}
\sum_{k \in \dT_n^i}\delta_k(\delta_{2k}\delta_{2k+1} -p^{(i)}(1,1)) =\sum_{k \in \dT_n}\underbrace{\ind{k \in \dT^i}\delta_k}_{H_k}\underbrace{\ind{k \in \dT^i}(\zeta^0_k\zeta^1_k-p^{(i)}(1,1))}_{G_k},
\end{equation*}
and Lemma~\ref{lemcvmart} again. \hspace{\stretch{1}}
$\Box$

\subsection{Laws of large numbers for the noise process}
\label{section noise}
We need to establish strong laws of large numbers for the noise sequence $(\varepsilon_n)$ restricted to the observed indices.

\begin{Lemma}\label{lemLLNeps}
Under assumptions \emph{(\textbf{HN.1})}, \emph{(\textbf{HN.2})}, \emph{(\textbf{HO})}, \emph{(\textbf{HI})} and for $i\in\{0,1\}$, one has
\begin{equation*}
\lim_{n \rightarrow + \infty}  \frac{1}{\pi^n}\sum_{k \in \dT_{n-1}}\delta_{2k+i}\veps_{2k+i}=0
\hspace{1cm}\text{a.s.}
\end{equation*}
\end{Lemma}

\noindent\textbf{Proof:} Set
\begin{equation*}
P_n^i=\sum_{k \in \dT_{n}}\underbrace{\delta_{2k+i}}_{H_k}\underbrace{\veps_{2k+i}}_{G_k}.
\end{equation*}
We use Lemma~\ref{lemcvmart}, with $\cG=(\dF^{\cO}_{n+1})$. Assumption (i)  is obvious. For $k \in \dG_{n+1}^i$, we have $\dE[G_k | \dF^{\cO}_{n+1}]=0$ and $\dE[G_k^2 | \dF^{\cO}_{n+1}]=\sigma^2$, and  $\dE[G_k G_p| \dF^{\cO}_{n+1}]=0$, for $k \neq p \in \dG_{n+1}^i$. Finally, we turn to assumption (iii):
\begin{equation*}
\sum_{k \in \dT_n} H_k^2 = \sum_{k \in \dT_n}\delta^2_{2k+i} = \sum_{\ell=1}^{n+1}Z_{\ell}^i= \cO(\pi^n),
\end{equation*}
as $n$ tends to infinity, thanks to  Eq.~\reff{limGW}. We obtain the result.\hspace{\stretch{1}}
\hspace{\stretch{1}}$\Box$\\

\begin{Lemma}\label{lemLLNeps2}
Under assumptions  \emph{(\textbf{HN.1})}, \emph{(\textbf{HN.2})}, \emph{(\textbf{HO})}, \emph{(\textbf{HI})} and
for $ i \in \{0,1\}$, one has
\begin{eqnarray*}
\lim_{n \rightarrow + \infty}\frac{1}{\pi^n}\sum_{k
\in\dT^i_{n}\backslash\dT_0}\veps_k^2\delta_k &= & \sigma^2\frac{\pi}{\pi-1} Wz^i
\quad \text{a.s.} \\
\lim_{n \rightarrow + \infty}\frac{1}{\pi^n}\sum_{k
\in\dT^i_{n}\backslash\dT_0}\!\!\!\delta_{2k}\delta_{2k+1}\veps_{2k}\veps_{2k+1}  &= &\rho p^{(i)}(1,1)\frac{ \pi }{\pi-1} Wz^i
\quad \text{a.s.} 
\end{eqnarray*}
\end{Lemma}

\noindent\textbf{Proof:} In order to prove the first convergence, we apply again Lemma \ref{lemcvmart} to the $\dF^{\cO}$-martingale:
\begin{equation*}
Q_n=\sum_{k \in\dT^i_{n}\backslash\dT_0}(\veps_k^2 - \sigma^2)\delta_k = \sum_{k \in\dT_{n}\backslash\dT_0}\underbrace{\ind{k \in \dT^i}\delta_{k}}_{H_k}\underbrace{\ind{k \in \dT^i}(\veps_k^2 - \sigma^2)}_{G_k},
\end{equation*}
Under (\textbf{HN.1}), (\textbf{HN.2}), we have $\dE[G_k | \dF^{\cO}_{n}]=0$ and $\dE[G_k^2 | \dF^{\cO}_{n}]=\tau^4-\sigma^4$, and  $\dE[G_k G_p| \dF^{\cO}_{n}]=0$, for $k \neq p \in \dG_n$. Thanks to  Eq.~\reff{limGW}, we have:
\begin{equation*}
\frac{1}{\pi^n} \sum_{k \in \dT^i_n}\delta_{k} =  \frac{1}{\pi^n} \sum_{\ell=1}^nZ_{\ell}^i \xrightarrow[n\rightarrow\infty]{} \frac{\pi}{\pi-1} Wz^i \quad \text{a.s.}
\end{equation*}which both implies assumption (iii) and  the first convergence.
To prove the second convergence, we write 
 \begin{eqnarray*}
\lefteqn{\frac{1}{\pi^n}\sum_{k\in\dT^i_{n}\backslash\dT_0}\!\!\!\delta_{2k}\delta_{2k+1}\veps_{2k}\veps_{2k+1}} \\
& =  &  \frac{1}{\pi^n} \sum_{k\in\dT_{n}\backslash\dT_0}\!\!\!\underbrace{\ind{k \in \dT^i}\delta_{2k}\delta_{2k+1}}_{H_k}\underbrace{\ind{k \in \dT^i}(\veps_{2k}\veps_{2k+1}-\rho)}_{G_k} 
  +  \frac{1}{\pi^n}\rho \sum_{k\in\dT^i_{n}\backslash\dT_0}\!\!\!\delta_{2k}\delta_{2k+1}
    \end{eqnarray*}
We use Lemma \ref{lemcvmart} to prove that the first term converges to $0$ ; Lemma \ref{lemSumdelta} gives the limit of the second term.\hspace{\stretch{1}}
$\Box$\\

\begin{Corollary}\label{CorLNeps2}
Under assumptions  \emph{(\textbf{HN.1})}, \emph{(\textbf{HN.2})}, \emph{(\textbf{HO})}, \emph{(\textbf{HI})} and for $ i \in \{0,1\}$, one has
\begin{eqnarray*}
\lim_{n \rightarrow + \infty}\frac{1}{\pi^n}\sum_{k \in\dT^i_{n}\backslash\dT_0}\veps_k^2\delta_{2k+j} & = & \sigma^2p_{ij}\frac{\pi}{\pi-1} Wz^i  \hspace{1cm}\text{a.s.}\\
\lim_{n \rightarrow + \infty} \frac{1}{\pi^n}\sum_{k \in \dT_{n}\backslash\dT_0}\delta_{2k}\delta_{2k+1}\veps_{2k}\veps_{2k+1} &  = & \rho  \bar{p}(1,1)\frac{\pi}{\pi-1} W 
\hspace{1cm}\text{a.s.}
\end{eqnarray*}
\end{Corollary}

\noindent\textbf{Proof:} The proof of the first limit is similar to the preceding ones, using the decomposition $\delta_{2k+j}=\delta_k\zeta_k^j$ and the properties of the sequence $(\zeta_n^j)$. Using Lemma~\ref{lemLLNeps2} the second one is straightforward. 
\hspace{\stretch{1}}$\Box$\\

\begin{Lemma}\label{lemLLEps4}
Under assumptions  \emph{(\textbf{HN.1})}, \emph{(\textbf{HN.2})}, \emph{(\textbf{HO})}, \emph{(\textbf{HI})} and for $i\in\{0,1\}$, we have
\begin{eqnarray*}
\lim_{n \rightarrow + \infty}\frac{1}{\pi^n}\sum_{k \in\dT_n^i
\backslash \dT_0}\delta_k\veps_k^4&=&\tau^4 \frac{\pi}{\pi-1} Wz^i\hspace{1cm}\text{a.s.}\\
\lim_{n \rightarrow +
\infty}\frac{1}{\pi^n}\sum_{k
\in\dT_{n-1}^i}\delta_{2k}\delta_{2k+1}\veps_{2k}^2\veps_{2k+1}^2&=&\nu^2 \tau^4 p^{(i)}(1,1)\frac{\pi}{\pi-1} Wz^i \hspace{1cm}\text{a.s.}
\end{eqnarray*}
\end{Lemma}

\noindent\textbf{Proof :} The proof 
follows essentially the same lines as the proof of
Lemma~\ref{lemLLNeps2} using the square integrable real martingales
\begin{equation*}
Q_n  =  \sum_{k \in\dT_n^i\backslash \dT_0}\delta_k(\veps_i^4-\tau^4),\quad\text{and}\quad
R_n  =  \sum_{k \in\dT_n^i \backslash \dT_0}
\delta_{2j}\delta_{2j+1}(\veps_{2j}^2\veps_{2j+1}^2-\nu^2 \tau^4).
\end{equation*}
It is therefore left to the reader.\hspace{\stretch{1}}$\Box$

\section{Convergence of the increasing process}
\label{section proof mainlemma}
We can now turn to the proof of our keystone result, the convergence of the increasing process of the main martingale $(\bs{M}_n)$.
\subsection{Preliminary results}
\label{prelim}
We first need an upper bound of the normalized sums of the $\delta_{2n+i}X_n^2$, and $\delta_{2n}\delta_{2n+1}X_n^2$ before being able to deduce their limits.

\begin{Lemma}\label{lem SumX}
Under assumptions \emph{(\textbf{HN.1})}, \emph{(\textbf{HN.2})}, \emph{(\textbf{HI})}  and \emph{(\textbf{HO})}, and for $i\in\{0,1\}$, we have
\begin{equation*}
\sum_{k \in \dT_n}\delta_{2k+i}X_k^2 = \cO (\pi^n) \quad \text{and} \quad
\sum_{k \in \dT_n}\delta_{2k}\delta_{2k+1}X_k^2 = \cO (\pi^n)\quad \text{a.s.}
\end{equation*}
\end{Lemma}

\noindent\textbf{Proof:} In all the sequel, for all $k\geq 1$,
define $a_{2k}=a$, $b_{2k}=b$, $a_{2k+1}=c$, $b_{2k+1}=d$ and
$\eta_k = a_k + \veps_k$ with the convention that $\eta_1=0$. It
follows from a recursive application of relation~(\ref{defbar}) that,
for all $k\geq 1$,
\begin{equation}\label{receq}
X_k = \Big(\prod_{\ell=0}^{r_k-1}b_{[\frac{k}{2^{\ell}}]}\Big)X_1 + \sum_{\ell=0}^{r_k-1}
\Big(\prod_{p=0}^{\ell-1}b_{[\frac{k}{2^p}]}\Big)\eta_{[\frac{k}{2^{\ell}}] },
\end{equation}
with the convention that an empty product equals $1$. Set $\alpha=\max(|a|, |c|) $, $\beta=\max(|b|, |d|)$ and notice that $0<\beta <1$.
The proof of Lemma~A.5 in \cite{BSG09} yields
\begin{eqnarray}
\sum_{k \in \dT_n \backslash \dT_0}\!\!\!\!\! \delta_{2k+i}X_k^2
& \leq & \frac{4}{1-\beta}\sum_{k \in \dT_n \backslash \dT_0}\!\!\!\!\!\delta_{2k+i}
\sum_{\ell=0}^{r_k-1}\beta^j\veps^2_{[\frac{k}{2^{\ell}}] }
+ \frac{4\alpha^2}{1-\beta}\sum_{k \in \dT_n \backslash \dT_0}\!\!\!\!\!\delta_{2k+i}\sum_{\ell=0}^{r_k-1}\beta^{\ell}\nonumber \\
&&+ 2X_1^2
\sum_{k \in \dT_n \backslash \dT_0}\delta_{2k+i}\beta^{2r_k}, \nonumber\\
& \leq &
\frac{4}{1-\beta}A_n^i+\frac{4\alpha^2}{1-\beta}B_n^i+2X_1^2C_n^i,
\label{decosumX}
\end{eqnarray}
where, for $i \in \{0,1\}$,
\begin{equation*}
A_n^i=\!\!\!\sum_{k \in \dT_n \backslash \dT_0}\!\!\!\!\!
\delta_{2k+i}\sum_{\ell=0}^{r_k-1}\beta^{\ell}\veps^2_{[\frac{k}{2^{\ell}}]},\hspace{0.3cm}
B_n^i=\!\!\!\sum_{k \in \dT_n \backslash
\dT_0}\!\!\!\!\!\delta_{2k+i}\sum_{\ell=0}^{r_k-1}\beta^{\ell},\hspace{0.3cm}
C_n^i=\!\!\!\sum_{k \in \dT_n \backslash
\dT_0}\!\!\!\!\!\delta_{2k+i}\beta^{2r_k}.
\end{equation*}
The last two terms above are readily evaluated by splitting the sums ge\-ne\-ra\-tion-wise.
Indeed, the last term can be rewritten as
\begin{equation*}
C_n^i=\sum_{\ell=1}^n\sum_{k \in \dG_{\ell}}\delta_{2k+i}\beta^{2\ell}= \sum_{\ell=1}^n \beta^{2\ell} Z^i_{\ell+1} = \pi^n \sum_{\ell=1}^n (\pi^{-1})^{n-\ell} \big(\beta^{2\ell} \frac{Z^i_{\ell+1}}{\pi^{\ell}}\big).
\end{equation*}
We now use Lemma \ref{lemtoepbar} with $\bs{A}_n=\pi^{-n}$ and $\bs{X}_n=\beta^{2n}{Z^i_{n+1}}{\pi^{-n}}$. On the one hand, the series of $(\pi^{-n})$ converges to $\pi/(\pi-1)$ as $\pi>1$ by assumption; on the other hand, $\beta^{2n}$ tends to $0$ as $n$ tends to infinity as $\beta<1$, and ${Z^i_{n}}{\pi^{-n}}$ converges a.s. to $Wz^i$ according to Eq.~\reff{limGW}, hence $\beta^{2n}{Z^i_{n+1}}{\pi^{-n}}$ tends to $0$ as $n$ tends to infinity. Lemma \ref{lemtoepbar} thus yields
\begin{equation*}
\lim_{n \rightarrow \infty} \sum_{\ell=1}^n (\pi^{-1})^{n-\ell} \big(\beta^{2\ell} \frac{Z^i_{\ell+1}}{\pi^{\ell}}\big)= 0 \quad \text{and} \quad C_n^i = o(\pi^n)\hspace{1cm} \text{a.s.}
\end{equation*}
We now turn to the term $B_n^i$:
\begin{equation*}
 B_n^i =\sum_{\ell=1}^n\sum_{k \in \dG_{\ell}}\delta_{2k+i}\frac{1-\beta^{\ell}}{1-\beta}
 \leq  \frac{1}{(1-\beta)}\sum_{\ell=1}^n\sum_{k \in \dG_{\ell}}\delta_{2k+i} 
\leq \frac{|\dT_{n+1}^*|}{(1-\beta)} =\ \cO(\pi^n),
\end{equation*}
due to Lemma \ref{lemcard}. It remains to control the first term $A_n^i$. Note that $\veps_k$ appears in $A_n^i$ as many times as it has descendants up to the $n$-th generation, and its multiplicative factor for its $p$-th generation descendant $k$ is $\beta^p\delta_{2k}$. This leads to
\begin{equation*}
A_n^i = \sum_{\ell=1}^n\sum_{k\in\dG_{\ell}}\veps_k^2 \sum_{p=0}^{n-\ell}\beta^p \sum_{m=0}^{2^p-1}\delta_{2(2^pk+m)+i}.
\end{equation*}
Now, note that $\sum_{m=0}^{2^p-1}\delta_{2(2^pk+m)+i}=\delta_k\sum_{m=0}^{2^p-1}\delta_{2(2^pk+m)+i}$ is the number of descendants of type $i$ of individual $k$ after $p+1$ generations. We denote it $Z^i_{p+1}(k)$, and split $A_n^i$ the following way:
\begin{equation} \label{loc}
A_n^i = \sum_{\ell=1}^n\sum_{k\in\dG_{\ell}} \sigma^2 \sum_{p=0}^{n-\ell}\beta^p \delta_k Z^i_{p+1}(k) +\sum_{\ell=1}^n\sum_{k\in\dG_{\ell}}(\veps_k^2 - \sigma^2) \sum_{p=0}^{n-\ell}\beta^p \delta_k Z^i_{p+1}(k).
\end{equation}
We first deal with the second term of the above sum.
\begin{eqnarray*}
 \sum_{\ell=1}^n\sum_{k\in\dG_{\ell}}(\veps_k^2 - \sigma^2) \sum_{p=0}^{n-\ell}\beta^p \delta_k Z^i_{p+1}(k)
&=& \sum_{p=0}^{n-1}\beta^p\sum_{\ell=1}^{n-p}\sum_{k\in\dG_{\ell}}(\veps_k^2 - \sigma^2) \delta_k Z^i_{p+1}(k) \\
&=& \sum_{p=0}^{n-1}\beta^p \sum_{\ell=1}^{n-p} Y_{\ell,p}^i,
\end{eqnarray*}
where $Y_{\ell,p}^i= \sum_{k\in\dG_{\ell}}(\veps_k^2 - \sigma^2) \delta_k Z^i_{p+1}(k)$. Tedious but straightforward computations lead to the following expression for the second order moment of $Y_{\ell,p}^i$, relying on assumptions (\textbf{HI}), (\textbf{HN.1}) and (\textbf{HN.2}). We also use the fact that, for $k \in \dG_{\ell}$, conditionally to $\{\delta_k=1\}$, $Z^i_{p+1}(k)$ follows the same law as $Z^i_{p+1}$, and is independent of any $Z^i_{p+1}(k')$, for $k'\neq k \in \dG_{\ell}$.
\begin{eqnarray*}
 \dE[(Y_{\ell,p}^i)^2]
& = &(\tau^4-\sigma^4) \dE[Z_{\ell}^0+Z_{\ell}^1]\dE[(Z^i_{p+1})^2] \\
& & \qquad  + (\nu^2 \tau ^4 -\sigma^4) \dE[Z^i_{p+1}]^2 \dE\Big[\sum_{k \in \dG_{\ell-1}} \delta_{2k}\delta_{2k+1}\Big] \\
& \le & (\tau^4-\sigma^4) \dE[Z_{\ell}^0+Z_{\ell}^1] \Big(\dE[(Z^i_{p+1})^2] + \dE[Z^i_{p+1}]^2\Big),
\end{eqnarray*}
since $\sum_{k \in \dG_{\ell-1}} \delta_{2k}\delta_{2k+1} \le \sum_{k \in \dG_{\ell-1}} (\delta_{2k}+\delta_{2k+1}) = Z_{\ell}^0+Z_{\ell}^1$. Now, using results on the moments of a two-type Galton-Watson process (see e.g. \cite{Har63}), we know that $\dE[(Z^i_{p+1})^2]= \cO(\pi^{2p})$. Recall Eq.~\reff{limGW} to obtain that $\dE[(Y_{\ell,p}^i)^2] = \cO(\pi^{\ell} \pi^{2p})$, which immediately entails that $|Y_{\ell,p}^i| = o(\pi^{\alpha \ell} \pi^{\gamma p})$ a.s., for any $\alpha >1/2$ and $\gamma > 1$. We thus one gets
\begin{equation*}
 \sum_{p=0}^{n-1}\beta^p \sum_{\ell=1}^{n-p} Y_{\ell,p}^i = \cO((\beta \pi^{\gamma})^n) = \cO(\pi^n) \hspace{1cm} \text{a.s.},
\end{equation*}
since we can choose $\gamma$ close enough to $1$ to get $\beta \pi^{\gamma} \le \pi$, as $\beta <1$. We have thus proved that the second term in the sum in \reff{loc} is $\cO(\pi^n)$, we now turn to the first one
\begin{eqnarray*}
\lefteqn{ \sum_{\ell=1}^n\sum_{k\in\dG_{\ell}} \sigma^2 \sum_{p=0}^{n-\ell}\beta^p \delta_k Z^i_{p+1}(k) }\\
& = & \sigma^2 \sum_{\ell=1}^n \sum_{p=0}^{n-\ell}\beta^p \sum_{k\in\dG_{\ell}} \delta_k Z^i_{p+1}(k) \ =\ \sigma^2 \sum_{\ell=1}^n \sum_{p=0}^{n-\ell}\beta^p Z^i_{\ell+p+1} \\
& = & \sigma^2 \sum_{p=0}^{n-1} \beta^p \sum_{\ell=1}^{n-p} Z^i_{\ell+p+1} \ \le \ \sigma^2 \sum_{p=0}^{n-1} \beta^p |\dT_{n+1}^*| = \cO(\pi^n) \hspace{1cm} \text{a.s.}
\end{eqnarray*}
Finally, $A_n^i = \cO(\pi^n)$, and the first result of the Lemma is proved. The second result follows immediately from the remark that the second sum in Lemma~\ref{lem SumX} is clearly smaller than the first one.\hspace{\stretch{1}}
$\Box$\\

\begin{Lemma}\label{lemSumX4}
Under assumptions \emph{(\textbf{HN.1})}, \emph{(\textbf{HN.2})}, \emph{(\textbf{HI})}  and \emph{(\textbf{HO})}, and for $i\in\{0,1\}$, we have
\begin{equation*}
\sum_{k \in \dT_n}\delta_{2k+i}X_k^4 = \cO (\pi^n) \quad \text{and} \quad
\sum_{k \in \dT_n}\delta_{2k}\delta_{2k+1}X_k^4 = \cO (\pi^n) \quad \text{a.s.}
\end{equation*}
\end{Lemma}

\noindent\textbf{Proof:} The proof mimics that of Lemma~\ref{lem SumX}. Instead of Equation~(\ref{decosumX}),
we have
\begin{equation*}
\sum_{k \in \dT_n \backslash \dT_0}  \delta_{2k+i}X_k^4  \leq
\frac{64}{(1-\beta)^3}A_n^i+\frac{64\alpha^4}{(1-\beta)^3}B_n^i+
8X_1^4C_n^i
\end{equation*}with, for $i$ in $\{0,1\}$
\begin{equation*}
A_n^i=\!\!\!\sum_{k \in \dT_n \backslash \dT_0}\!\!\!\!\!
\delta_{2k+i}\sum_{\ell=0}^{r_k-1}\beta^j\veps^4_{[\frac{k}{2^{\ell}}]},\hspace{0.3cm}
B_n^i=\!\!\!\sum_{k \in \dT_n \backslash
\dT_0}\!\!\!\!\!\delta_{2k+i}\sum_{\ell=0}^{r_k-1}\beta^{\ell},\hspace{0.3cm}
C_n^i=\!\!\!\sum_{k \in \dT_n \backslash
\dT_0}\!\!\!\!\!\delta_{2k+i}\beta^{4r_k}.
\end{equation*}
We can easily prove that $(B_n^i+C_n^i)=\cO(\pi^n)$.
Therefore, we only need a sharper estimate for $A_n^i$. Via the same
lines as in the proof of Lemma~\ref{lem SumX}, but dealing with $\veps_k^4$ instead of $\veps_k^2$, we can show that
$A_n^i=\cO(\pi^n)$ a.s. which immediately yields the first result. The second one is obtained by remarking that the second sum is less than the first one.\hspace{\stretch{1}}$\Box$

\subsection{Asymptotic behavior of the sum of observed data}
\label{LGN X}
We now turn to the asymptotic behavior of the sums of the observed data. More precisely, set $H_n^i =\sum_{k\in\dT_n}\delta_{2k+i}X_k$, for $i$ in $\{0,1\}$, and $\bs{H}_n=(H_n^0,H_n^1)^t$. The following result gives the asymptotic behavior of $(\bs{H}_n)$.

\begin{Proposition}\label{propSumX}
Under assumptions \emph{(\textbf{HN.1})}, \emph{(\textbf{HN.2})}, \emph{(\textbf{HI})}  and \emph{(\textbf{HO})}, we have the convergence:
\begin{eqnarray*}
\lim_{n\rightarrow \infty} \frac{\bs{H}_{n}}{\pi^n}= \frac{\pi}{\pi-1} W  \bs{h} \hspace{1cm}\text{a.s.},
\end{eqnarray*}
where
\begin{equation*}
\bs{h}=\left(\begin{array}{c}
h^0\\
h^1
\end{array}\right)=(\rI_2 - \widetilde{\bs{P}}_1)^{-1}\bs{P}^t\left(\begin{array}{c}
az^0\\cz^1
\end{array}\right) \quad \text{and} \quad
\widetilde{\bs{P}}_1  =  \frac{1}{\pi}{\bs{P}^t}\left(\begin{array}{cc}
b&0\\0&d
\end{array}\right).
\end{equation*}
\end{Proposition}

\noindent\textbf{Proof:} We first prove that the sequence $(\bs{H}_n)$ satisfies a recursive property using
Equation~(\ref{defbar}).
\begin{eqnarray*}
H_n^0 & = & X_1\delta_2  +\!\!\! \sum_{k \in \dT_n^0}\!\!\left(a + bX_{[\frac{k}{2}] } + \veps_k\right)\delta_{2k}
+\!\!\! \sum_{k \in \dT_n^1\backslash \dT_0}\!\!\left(c + d X_{[\frac{k}{2}] } + \veps_k\right)\delta_{2k} \\
& = & X_1\delta_2 + a\!\!\sum_{k \in \dT_n^0}\!\!\delta_{2k} +  b\!\!\sum_{k \in \dT_n^0}\!\!X_{[\frac{k}{2}] }\delta_{2k}
+ c\!\!\sum_{k \in \dT_n^1\backslash \dT_0}\!\!\delta_{2k} + d\!\!\sum_{k \in \dT_n^1\backslash \dT_0}\!\!X_{[\frac{k}{2}] }\delta_{2k}\\
&& + \sum_{k \in \dT_n\backslash \dT_0}\veps_{k}\delta_{2k} \\
& = &  bp_{00} H_{n-1}^0 + dp_{10} H_{n-1}^1 + B^0_n,
\end{eqnarray*}
with
\begin{eqnarray*}
B^0_n & = & X_1\delta_2 + a\sum_{k \in \dT_n^0}\delta_{2k}
   + c\sum_{k \in \dT_n^1\backslash \dT_0}\delta_{2k} +\sum_{k \in \dT_n\backslash \dT_0}\veps_{k}\delta_{2k} \\
    &  & \quad
   + b\sum_{k \in \dT_{n-1}}X_k\delta_{2k} (\delta_{4k} - p_{00})
   +  d\sum_{k \in \dT_{n-1}}X_k \delta_{2k+1}(\delta_{4k+2} -  p_{10}).
\end{eqnarray*}
Similarly, we have
\begin{equation*}
H_n^1 = bp_{01} H_{n-1}^0 + dp_{11} H_{n-1}^1 + B^1_n,
\end{equation*}
with
\begin{eqnarray*}
B^1_n & = & X_1\delta_3 + a\sum_{k \in \dT_n^0}\delta_{2k+1} + c\sum_{k \in \dT_n^1\backslash \dT_0}\delta_{2k+1} +\sum_{k \in \dT_n\backslash \dT_0}\veps_{k}\delta_{2k+1} \\
&& + b\sum_{k \in \dT_{n-1}}X_k\delta_{2k} (\delta_{4k+1} - p_{01}) + d\sum_{k \in \dT_{n-1}}X_k \delta_{2k+1}(\delta_{4k+3} -  p_{11}).
\end{eqnarray*}
Let us denote  $\bs{B}_n=(B_n^0,B_n^1)^t$. The last equations yield in the matrix form:
\begin{eqnarray*}\label{recursH}
\frac{\bs{H}_n}{\pi^{n}} & = & \widetilde{\bs{P}}_1\frac{\bs{H}_{n-1}}{\pi^{n-1}} + \frac{\bs{B}_n}{\pi^{n}}\ = \
 \widetilde{\bs{P}}_1^{n}\bs{H}_0 + \sum_{k=1}^n\widetilde{\bs{P}}_1^{n-k}\frac{\bs{B}_k}{\pi^{k}},
\end{eqnarray*}
with
\begin{equation*}
\widetilde{\bs{P}}_1=\frac{1}{\pi}\left(\begin{array}{cc}
bp_{00}&dp_{10}\\
bp_{01}&dp_{11}
\end{array}\right)=\frac{1}{\pi}{\bs{P}^t}\left(\begin{array}{cc}
b&0\\0&d
\end{array}\right).
\end{equation*}
One has $\|\widetilde{\bs{P}}_1^n\|\leq \pi^{-n}\beta^n\|\bs{P}^n\|$. It is well known that $\pi^{-n}\bs{P}^n$ converges to a fixed matrix (see e.g. \cite{HJ90}) as $\bs{P}$ is a positive matrix with dominant eigenvalue $\pi$. Since $\beta<1$, the sequence $\widetilde{\bs{P}}_1^n$ thus converges to $0$ as $n$ tends to infinity. In addition, $\sum\|\widetilde{\bs{P}}_1^n\|$ is bounded, $\rI_2 - \widetilde{\bs{P}}_1$ is invertible and $\sum_{n \ge 0} \widetilde{\bs{P}}_1^n$ converges to $(\rI_2 - \widetilde{\bs{P}}_1)^{-1}$. In order to use Lemma~\ref{lemtoepbar}, we need to compute the limit of $\bs{B}_n/\pi^n$. First, we prove that
\begin{equation} \label{cvsumepsdelt}
\sum_{k \in \dT_n\backslash \dT_0}\veps_{k}\delta_{2k+i}=o(\pi^n),
\end{equation}
for $i \in \{0,1\}$, thanks to Lemma~\ref{lemcvmart}. Indeed, set $\mathcal{G}=\dF^{\cO}$, $H_k=\delta_{2k+i}$, $G_k=\veps_k$. Thus hypothesis (i) of Lemma~\ref{lemcvmart} is obvious, (ii) comes from \textbf{(HN.1)} and \textbf{(HN.2)}. Finally, the last assumption (iii) holds, since
\begin{equation*}
 \sum_{k \in \dT_n\backslash \dT_0} \delta_{2k+i}^2 = \sum_{\ell=1}^{n+1} Z_{\ell}^i = \cO(\pi^n),
\end{equation*}
the last equality coming from \reff{limGW}, which holds thanks to \textbf{(HO)}. Now, we turn to the terms
$$\sum_{k \in \dT_n}X_k\delta_{2k+i}(\delta_{2(2k+i)+j} - p_{ij})= \sum_{k \in \dT_n}X_k\delta_{2k+i}(\zeta_{2k+i}^j - p_{ij}),$$
for $(i,j) \in \{0,1\}^2$. We use again Lemma~\ref{lemcvmart}, with the following setting: $(\cG_n)=(\cZ_{n+1} \vee \cF_{n+1})$, $H_k=X_k \delta_{2k+i}$, $G_k=\zeta_{2k+i}^j - p_{ij}$. For $k \in \dG_n$, we check that $X_k \delta_{2k+i}$ is $\cG_{n-1}$-measurable, since $X_k$ is $\cF_n$-measurable and $\delta_{2k+i}$ is $\cZ_n$-measurable. Next, because of \textbf{(HI)} and of the independence of the sequence $(\bs{\zeta}_k)$, $\dE[\zeta_{2k+i}^j - p_{ij} | \cZ_n \vee \cF_n]= 0$. The same independence hypothesis yields that $\dE[G_kG_p |\cZ_n \vee \cF_n] \neq 0$ only if $k=p$, and then equals $\sigma_{ij}^2$. Finally, 
\begin{equation*}
 \sum_{k \in \dT_n}(X_k\delta_{2k+i})^2 = \sum_{k \in \dT_n}X_k^2\delta_{2k+i} = \cO(\pi^n),
\end{equation*}
thanks to Lemma~\ref{lem SumX}. Now, Lemma~\ref{lemcvmart} allows to conclude that
\begin{equation}\label{cvint}
 \sum_{k \in \dT_n}X_k\delta_{2k+i}(\delta_{2(2k+i)+j} - p_{ij})= o(\pi^n),
\end{equation}
for $(i,j) \in \{0,1\}^2$. Next, Lemma~\ref{lemSumdelta} gives the limit of the term $\sum_{k \in \dT_n^i} \delta_{2k+j}$, for $(i,j) \in \{0,1\}^2$, so that we finally obtain:
\begin{equation*}
\lim_{n\rightarrow \infty} \frac{\bs{B}_{n}}{\pi^n} =  W\frac{\pi}{\pi-1} \left(\begin{array}{c}
a  z^0p_{00} + cz^1 p_{10}\\
a  z^0p_{01} + cz^1 p_{11}
\end{array}\right)
= W\frac{\pi}{\pi-1}\bs{P}^t\left(\begin{array}{c}az^0\\cz^1\end{array}\right) \hspace{1cm}\text{a.s.}
\end{equation*}
and we use Lemma~\ref{lemtoepbar} to conclude.\hspace{\stretch{1}}$ \Box$\\

\begin{Remark}\label{remSumX}
Putting together Proposition~\ref{propSumX} and Eq.~\reff{cvint} above, we immediately get that under the same assumptions as that of Proposition~\ref{propSumX}, 
\begin{equation*}
\lim_{n\rightarrow \infty}\frac{1}{\pi^n}\sum_{k \in\dT_n}X_k\delta_{2k+i}\delta_{2(2k+i)+j}
= \frac{\pi}{\pi-1}  h^i p_{ij}W \hspace{1cm}\text{a.s.}
\end{equation*}
for all $(i,j)\in\{0,1\}^2$, result we will use for the study of the limit of $\sum X_k^2 \delta_{2k+i}$. 
\end{Remark}

\subsection{Asymptotic behavior of the sum of squared observed data}
\label{carres}
We now turn to the asymptotic behavior of the sums of the squared observed data. Set $K_n^i =\sum_{k\in\dT_n}\delta_{2k+i}X^2_k$, for $i$ in $\{0,1\}$, and $\bs{K}_n=(K_n^0,K_n^1)^t$.  The following result gives the asymptotic behavior of $(\bs{K}_n)$.
\begin{Proposition}\label{propSumX2}
Under assumptions \emph{(\textbf{HN.1})}, \emph{(\textbf{HN.2})}, \emph{(\textbf{HI})}  and \emph{(\textbf{HO})}, we have the convergence:
\begin{eqnarray*}
\lim_{n\rightarrow \infty} \frac{\bs{K}_{n}}{\pi^n}= \frac{\pi}{\pi-1} W  \bs{k} \hspace{1cm}\text{a.s.},
\end{eqnarray*}
where
\begin{equation*}
\bs{k}= \left(\begin{array}{c}
k^0\\
k^1
\end{array}\right) =(\rI_2 - \widetilde{\bs{P}}_2)^{-1} \bs{P}^t \left(\begin{array}{c}
(a^2+\sigma^2)z^0+\frac{2}{\pi}ab h^0\\
(c^2+\sigma^2)z^1+\frac{2}{\pi}cd h^1
\end{array}\right),
\end{equation*}
and
\begin{equation*}
\widetilde{\bs{P}}_2  =  \frac{1}{\pi}{\bs{P}^t}\left(\begin{array}{cc}
b^2&0\\0&d^2
\end{array}\right).
\end{equation*}
\end{Proposition}

\noindent\textbf{Proof:} We use again Equation~(\ref{defbar}) to prove a recursive property for the sequence $(\bs{K}_n)$.
Following the same lines as in the proof of Proposition~\ref{propSumX}, we obtain:
\begin{eqnarray*}\label{recursK}
\frac{\bs{K}_n}{\pi^{n}} & = & \widetilde{\bs{P}}_2\frac{\bs{K}_{n-1}}{\pi^{n-1}} + \frac{\bs{C}_n}{\pi^{n}}\ = \ \widetilde{\bs{P}}_2^{n}\bs{K}_0 + \sum_{\ell=1}^n\widetilde{\bs{P}}_2^{n-\ell}\frac{\bs{C}_{\ell}}{\pi^{\ell}},
\end{eqnarray*}
where $\bs{C}_n=(C_n^0,C_n^1)^t$ is defined by
\begin{eqnarray*}
C^i_n & = & X^2_1\delta_{2+i} + a^2\sum_{k \in \dT_n^0}\delta_{2k+i}+ b^2\sum_{k \in \dT_{n-1}}X^2_k\delta_{2k} (\delta_{4k+i}
 - p_{0i})\\
&& + 2ab\!\!\!\sum_{k \in \dT_{n-1}}\!\!\!X_k\delta_{2k}\delta_{4k+i}
+2a\!\sum_{k \in \dT_n^0}\!\veps_k\delta_{2k+i}
+2b\!\sum_{k \in \dT_n^0}\!X_{[\frac{k}{2}] }\veps_k\delta_{2k+i}\\
&&+\!\!\!\!\sum_{k \in \dT_n\backslash \dT_0}\!\!\!\veps^2_{k}\delta_{2k+i} + c^2\!\!\!\!\sum_{k \in \dT_n^1\backslash \dT_0}\!\!\!\delta_{2k+i}+d^2\!\!\!\!\sum_{k \in \dT_{n-1}}\!\!\!X^2_k\delta_{2k+1}(\delta_{4k+2+i} -  p_{1i})  \\
&& + 2cd\!\!\!\sum_{k \in \dT_{n-1}}\!\!\!X_k\delta_{2k+1}\delta_{4k+2+i}  +2c\!\!\!\!\sum_{k \in \dT_n^1\backslash \dT_0}\!\!\!\veps_k\delta_{2k+i}
+2d\!\!\sum_{k \in \dT_n^1\backslash \dT_0}\!X_{[\frac{k}{2}] }\veps_k\delta_{2k+i},
\end{eqnarray*}
for $i \in \{0,1\}$.
Note that $\|\widetilde{\bs{P}}_2^n\|\leq \pi^{-n}\beta^{2n}\|\bs{P}^n\|$, so that $\widetilde{\bs{P}}_2^n$ converges to $0$. In addition, $\sum\|\widetilde{\bs{P}}_2^n\|$ is bounded, $\rI_2 - \widetilde{\bs{P}}_2$ is invertible and $\sum_{n \ge 0}\widetilde{\bs{P}}_2^n$ converges to $(\rI_2 - \widetilde{\bs{P}}_2)^{-1}$. In order to use Lemma~\ref{lemtoepbar}, we have to compute the limit of $\bs{C}_n/\pi^n$. Following the proof of \reff{cvsumepsdelt}, we already have, for $(i,j) \in \{0,1\}^2$,
\begin{equation*} 
\sum_{k \in \dT_n^j} \veps_k\delta_{2k+i} = o(\pi^n) \quad \text{a.s.}
\end{equation*}
We now turn to the terms $\sum_{k \in \dT_{n-1}}X^2_k\delta_{2k+i} (\delta_{2(2k+i)+j} - p_{ij})$,for $(i,j) \in \{0,1\}^2$. To deal with these terms, we use Lemma~\ref{lemcvmart} with the same setting we used to prove Eq.~\reff{cvint}, except that we replace $X_k$ with $X_k^2$. Assumptions (i) and (ii) of Lemma~\ref{lemcvmart} have thus already been checked, and regarding (iii), we have $\sum_{k \in \dT_{n-1}} X_k^4 \delta_{2k+i} = \cO(\pi^n)$ a.s. thanks to Lemma~\ref{lemSumX4}. We conclude that
\begin{equation*} 
\sum_{k \in \dT_{n-1}}X^2_k\delta_{2k+i} (\delta_{2(2k+i)+j}
 - p_{ij}) = o(\pi^n) \quad \text{a.s.}
\end{equation*}
Next, we study $\sum_{k \in \dT_n^i} X_{[\frac{k}{2}]} \veps_k \delta_{2k+j}$, for $(i,j) \in \{0,1\}^2$. We use the same martingale tool, so to speak Lemma~\ref{lemcvmart}, with $\mathcal{G}= \dF^{\mathcal O}$, $H_k = X_{[\frac{k}{2}]} \delta_{2k+j} \ind{k \in \dT^i}$ and $G_k = \veps_k$. Assumptions (i) and (ii) are easily checked, and since 
$$\sum_{k \in \dT_n^i} X_{[\frac{k}{2}]}^2 \delta_{2k+j} = \sum_{k \in \dT_{n-1}} X_k^2 \delta_{2(2k+i)+j} \le\sum_{k \in \dT_{n-1}} X_k^2 \delta_{2k+i} = \cO(\pi^n),$$
the last equality coming from Lemma~\ref{lem SumX}, assumption (iii) is satisfied and
\begin{equation*} 
 \sum_{k \in \dT_n^i} X_{[\frac{k}{2}]} \veps_k \delta_{2k+j} = o(\pi^n) \quad \text{a.s.}
\end{equation*}
Now, Corollary~\ref{CorLNeps2} yields that for $i \in \{0,1\}$,
\begin{equation*}
\lim_{n \rightarrow \infty} \frac{1}{\pi^n}\sum_{k \in \dT_n\backslash \dT_0} \veps^2_{k}\delta_{2k+i} = \sigma^2 (p_{0i}z^0+p_{1i}z^1) \frac{\pi}{\pi-1}W \quad \text{a.s.}
\end{equation*}
Finally, Remark~\ref{remSumX} gives the limit of $\pi^{-n}\sum_{k \in \dT_{n-1}} X_k\delta_{2k+i}\delta_{2(2k+i)+j}$, and Lemma~\ref{lemSumdelta} that of $\pi^{-n}\sum_{k \in \dT_n^j}\delta_{2k+i}$, so that we finally obtain
\begin{equation*}
\lim_{n\rightarrow \infty} \frac{\bs{C}_{n}}{\pi^n}= \frac{W\pi}{\pi-1}
\left(\begin{array}{cc}
         p_{00} & p_{10}\\
         p_{01} & p_{11}
        \end{array}\right)
\times 
\left(\begin{array}{c}
       (a^2+\sigma^2)  z^0 + \frac{2}{\pi}ab h^0 \\
       (c^2+\sigma^2)z^1 + \frac{2}{\pi}cd h^1
      \end{array}\right) \text{a.s.}
\end{equation*}
And we conclude using Lemma~\ref{lemtoepbar} again.\hspace{\stretch{1}}$ \Box$\\

Propositions \ref{propSumX} and \ref{propSumX2} together with Equation~(\ref{limGW}) give the asymptotic behavior of the matrices $\bs{S}_n^0$ and $\bs{S}_n^1$. The next result gives the behavior of matrix $\bs{S}_n^{0,1}$ given through the quantities $\sum_{k\in\dT_n}\delta_{2k}\delta_{2k+1}X_k$ and $\sum_{k\in\dT_n}\delta_{2k}\delta_{2k+1}X^2_k$. It is an easy consequence of Propositions \ref{propSumX} and \ref{propSumX2}, together with Lemma~\ref{lemSumdelta} for the first limit.

\subsection{Asymptotic behavior of covariance terms}
\label{covariance}
Finally, we turn to the asymptotic behavior of the covariance terms, which are involved in matrix $\bs{S}_n^{0,1}$. We thus define $H_n^{01} = \sum_{k\in\dT_n}\delta_{2k}\delta_{2k+1}X_k$ and $K_n^{01} = \sum_{k\in\dT_n}\delta_{2k}\delta_{2k+1}X_k^2$.

\begin{Proposition}\label{propS 01}
Under assumptions \emph{(\textbf{HN.1})}, \emph{(\textbf{HN.2})}, \emph{(\textbf{HO})} and \emph{(\textbf{HI})}, we have the almost sure convergences:
\begin{equation*}
\lim_{n\rightarrow \infty} \frac{1}{\pi^n} \sum_{k\in\dT_n} \delta_{2k}\delta_{2k+1}
=\frac{\pi}{\pi-1} W  \bar{p}(1,1),
\end{equation*}
\begin{equation*}
\lim_{n\rightarrow \infty} \frac{H_n^{01}}{\pi^n} = \frac{\pi}{\pi -1} W h^{0,1}
\quad \text{and} \quad
\lim_{n\rightarrow \infty} \frac{K_n^{01}}{\pi^n} = \frac{\pi}{\pi -1} W  k^{0,1},
\end{equation*}
where 
\begin{eqnarray}
\bar{p}(1,1) & = &  p^{(0)}(1,1) z^0 + p^{(1)}(1,1) z^1,\label{def pbar}\\
h^{0,1} & = & p^{(0)}(1,1) \left(a z^0 + b\frac{h^0}{\pi} \right) + p^{(1)}(1,1)\left(c z^1 +d \frac{h^1}{\pi}\right), \nonumber\\
k^{0,1} & = &  p^{(0)}(1,1) \left(a^2 z^0 + b^2 \frac{k^0}{\pi} + 2ab \frac{h^0}{\pi}\right)\nonumber\\ 
&   & + p^{(1)}(1,1) \left(c^2 z^1 + d^2 \frac{k^1}{\pi} + 2cd \frac{h^1}{\pi}\right) + \sigma^2 \bar{p}(1,1).\nonumber
\end{eqnarray}
\end{Proposition}

\noindent\textbf{Proof:} The first limit is a consequence of Lemma \ref{lemSumdelta}. Next, using Eq.~(\ref{defbar}) we obtain  ${H_n^{01}}{\pi^{-n}} $ and  ${K_n^{01}}{\pi^{-n}} $ in terms  of ${\pi^{-n}} \sum_{k\in\dT^i_{n-1}} \delta_{k}$, ${H_{n-1}^{i}}{\pi^{-n}} $ and ${K_{n-1}^{i}}{\pi^{-n}} $  and the result follows from Propositions \ref{propSumX} and \ref{propSumX2}.\hspace{\stretch{1}}$ \Box$\\

\noindent\textbf{Proof of Proposition~\ref{mainlemma}:} We are now in a position to complete the proof of  Proposition~\ref{mainlemma}. Simply notice that we have proved in Propositions~\ref{propSumX}, \ref{propSumX2} and \ref{propS 01} all the wished convergences, except that we normalized the sums with $\pi^n$. Thanks to Lemma~\ref{lemcard}, we end the proof.\hspace{\stretch{1}}$ \Box$

\begin{Remark} In the case of fully observed date, the matrix $\bs{P}$ is a $2\times2$ matrix with all entries equal to $1$, $\pi$ equals $2$ and the normalized eigenvector $\bs{z}$ equals $(1/2,1/2)$. One can check that in that case, our limits correspond to those of~\cite{BSG09}. 
\end{Remark}

\section{Asymptotic behavior of the main martingale}
\label{section proof M}
Theorem \ref{th Mn} is a strong law of large numbers for the martingale $(\bs{M}_n)$.
The standard strong law for martingales is unhelpful here. Indeed, it is valid for martingales that can be decomposed in a sum of the form $\sum_{\ell=1}^n\bs{\Psi}_{\ell-1}\bs{\xi}_{\ell}$ where $(\bs{\Psi}_{\ell})$ is predictable and $(\bs{\xi}_{\ell})$ is a martingale difference sequence. In addition, $(\bs{\Psi}_{\ell})$ and $(\bs{\xi}_{\ell})$ are required to be sequences of \emph{fixed-size} vectors. Such a decomposition with fixed-sized vectors is impossible in our context (see Lemma~\ref{lemlimB}), essentially because the number of observed data in each generation asymptotically grows exponentially fast as $\pi^n$.
Consequently, we are led to propose a new strong law of large numbers for $(\bs{M}_n)$, adapted to our framework.

\smallskip

For all $n\geq 1$, let $\cV_{n}=\bs{M}_n^t\bs{\Sigma}_{n-1}^{-1}\bs{M}_n$ where $\bs{\Sigma}_n$ is defined in Section~\ref{section def LS}. First of all, we have
\begin{eqnarray*}
\lefteqn{\cV_{n+1}}\\
&\!=\!&(\bs{M}_n+\Delta \bs{M}_{n+1})^t\bs{\Sigma}_n^{-1}(\bs{M}_n+\Delta \bs{M}_{n+1}),\\
&\!=\!&\cV_n\!-\!\bs{M}_n^{t}(\bs{\Sigma}_{n-1}^{-1}\!-\!\bs{\Sigma}_n^{-1})\bs{M}_n\!+\!
2\bs{M}_n^t\bs{\Sigma}_n^{-1}\Delta \bs{M}_{n+1}\!+\!\Delta \bs{M}_{n+1}^t\bs{\Sigma}_n^{-1}\Delta \bs{M}_{n+1}.
\end{eqnarray*}
Note that $\bs{M}_n^t\bs{\Sigma}_n^{-1}\Delta \bs{M}_{n+1}$ and $\Delta\bs{M}_n^t\bs{\Sigma}_n^{-1} \bs{M}_{n+1}$ are scalars, hence they are equal to their own transpose and as a result, one has $\bs{M}_n^t\bs{\Sigma}_n^{-1}\Delta \bs{M}_{n+1}=\Delta\bs{M}_n^t\bs{\Sigma}_n^{-1} \bs{M}_{n+1}$.
By summing over the identity above, we obtain the main decomposition
\begin{equation}\label{maindecomart}
\cV_{n+1}+\cA_n=\cV_1+\cB_{n+1}+\cW_{n+1},
\end{equation}
where
\begin{eqnarray*}
\cA_n=\sum_{\ell=1}^n \bs{M}_{\ell}^{t}(\bs{\Sigma}_{\ell-1}^{-1}-\bs{\Sigma}_{\ell}^{-1})\bs{M}_{\ell},
\end{eqnarray*}
\vspace{-0.8cm}
\begin{eqnarray*}
\cB_{n+1}=2\sum_{\ell=1}^n \bs{M}_{\ell}^t\bs{\Sigma}_{\ell}^{-1}\Delta \bs{M}_{\ell+1},\quad
\cW_{n+1}=\sum_{\ell=1}^n \Delta \bs{M}_{\ell+1}^t\bs{\Sigma}_{\ell}^{-1}\Delta \bs{M}_{\ell+1}.
\end{eqnarray*}
The asymptotic behavior of the left-hand side of (\ref{maindecomart}) is as follows.
\begin{Proposition}\label{lem lim V+A}
Under assumptions \emph{(\textbf{HN.1})}, \emph{(\textbf{HN.2})}, \emph{(\textbf{HO})} and \emph{(\textbf{HI})}, we have
\begin{equation*}
\lim_{n \rightarrow + \infty}\ind{|\dG_n^*|>0}\frac{\cV_{n+1}+\cA_n}{n} = \frac{4(\pi-1)}{\pi}\sigma^2\indnex
\hspace{1cm}\text{a.s.}
\end{equation*}
\end{Proposition}

\noindent\textbf{Proof :} Thanks to the laws of large numbers derived in Sections~\ref{LLN} and \ref{section proof mainlemma},  the proof of Proposition~\ref{lem lim V+A} follows essentially the same lines as \cite{BSG09} and is given in Appendix~\ref{appendixB}.
\hspace{\stretch{1}}$\Box$\\

Since $(\cV_{n})$ and $(\cA_n)$ are two sequences of non negative real numbers, 
Proposition~\ref{lem lim V+A} yields that $\ind{|\dG_n^*|>0}\cV_{n}=\cO(n)$ a.s. which proves Equation~(\ref{th Mn1}). We now turn to the proof of Equation~(\ref{th Mn3}). We start with a sharp rate of convergence for $(\bs{M}_n)$.

\begin{Proposition}\label{lemlimM}
Under assumptions \emph{(\textbf{HN.1})}, \emph{(\textbf{HN.2})}, \emph{(\textbf{HO})} and \emph{(\textbf{HI})}, we, we have, for all $\eta>1/2$,
\begin{equation*}
\ind{|\dG_n^*|>0}\parallel \bs{M}_n \parallel^2=o(|\dT_{n-1}^*| n^\eta)
\hspace{1cm}\text{a.s.}
\end{equation*}
\end{Proposition}

\noindent\textbf{Proof :} The result is obvious on $\cE$. On $\overline{\cE}$, the proof follows again the same lines as \cite{BSG09} thanks to the laws of large numbers derived in Sections~\ref{LLN} and \ref{section proof mainlemma}. It is given in Appendix~\ref{appendixC}. \hspace{\stretch{1}}$\Box$\\

A direct application of Proposition~\ref{lemlimM} ensures that $\ind{|\dG_n^*|>0}\cV_n=o(n^{\eta})$ a.s.
for all $\eta>1/2$. Hence, Proposition~\ref{lem lim V+A} immediately leads to the following result.

\begin{Corollary}\label{corlimA}
Under assumptions \emph{(\textbf{HN.1})}, \emph{(\textbf{HN.2})}, \emph{(\textbf{HO})} and \emph{(\textbf{HI})},  we have
\begin{equation*}
\lim_{n \rightarrow + \infty}\ind{|\dG_n^*|>0}\frac{\cA_n}{n} = \frac{4(\pi-1)}{\pi}\sigma^2\indnex
\hspace{1cm}\text{a.s.}
\end{equation*}
\end{Corollary}

\noindent\textbf{Proof of Result~(\ref{th Mn3}) of Theorem \ref{th Mn}:} First of all, $\cA_n$ may be rewritten as
\begin{equation*}
\cA_n=\sum_{\ell=1}^n \bs{M}_{\ell}^{t}(\bs{\Sigma}_{\ell-1}^{-1}-\bs{\Sigma}_{\ell}^{-1})\bs{M}_{\ell}
= \sum_{\ell=1}^n\bs{M}_{\ell}^{t}\bs{\Sigma}_{\ell-1}^{-1/2}\bs{\Delta}_{\ell} \bs{\Sigma}_{\ell-1}^{-1/2}\bs{M}_{\ell}
\end{equation*}
where $\bs{\Delta}_n=\rI_{4} - \bs{\Sigma}_{n-1}^{1/2}\bs{\Sigma}_n^{-1} \bs{\Sigma}_{n-1}^{1/2}$.
Thanks to Corollary~\ref{cv Gamman}, we know that
\begin{equation*}
\lim_{n\rightarrow\infty}\ind{|\dG_n^*|>0}\bs{\Delta}_n=\frac{\pi-1}{\pi}\rI_{4}\indnex
\hspace{1cm}\text{a.s.}
\end{equation*}
Besides, Corollary~\ref{corlimA} yields that $\cA_n \sim n \frac{\pi-1}{\pi} 4 \sigma^2$ a.s. on $\overline{\cE}$. Plugging these two results into the equality
\begin{equation*} 
\cA_n=\frac{\pi-1}{\pi}\sum_{\ell=1}^n\bs{M}_{\ell}^{t}\bs{\Sigma}_{\ell-1}^{-1} \bs{M}_{\ell} 
+ \sum_{\ell=1}^n\bs{M}_{\ell}^{t}\bs{\Sigma}_{\ell-1}^{-1/2}(\bs{\Delta}_{\ell}-\frac{\pi-1}{\pi}\rI_4)\bs{\Sigma}_{\ell-1}^{-1/2}\bs{M}_{\ell}
\end{equation*}
gives that $\sum_{\ell=1}^n\bs{M}_{\ell}^{t}\bs{\Sigma}_{\ell-1}^{-1} \bs{M}_{\ell} \sim \cA_n \frac{\pi}{\pi-1}$ a.s. on $\overline{\cE}$
and convergence~(\ref{th Mn3}) directly follows.\hspace{\stretch{1}}$
\Box$

\section{Proof of the main results}
\label{section main proof}
We can now proceed to proving our main results.
\subsection{Strong consistency for $\wh{\bs{\theta}}_n$}
\label{section main proof1}
Theorem \ref{thmaptheta} is a direct consequence of Theorem~\ref{th Mn}.\\

\noindent\textbf{Proof of result~(\ref{thmaptheta1}) of Theorem~\ref{thmaptheta}:} Recall that $\cV_{n}=\bs{M}_n^t\bs{\Sigma}_{n-1}^{-1}\bs{M}_n$. It clearly follows from Equation~(\ref{thetadiff}) that
\begin{equation*}
\cV_n=(\wh{\bs{\theta}}_n-\bs{\theta})^t\bs{\Sigma}_{n-1}(\wh{\bs{\theta}}_n-\bs{\theta)}.
\end{equation*}
Consequently, the asymptotic behavior of $\bs{\wh{\theta}}_n-\bs{\theta}$ is clearly related to the one of $\cV_n$. More precisely, we can deduce from Corollary~\ref{cv Gamman} and the fact that the eigenvalues of a matrix are continuous functions of its coefficients the following result
\begin{equation*}
\lim_{n\rightarrow\infty}\ind{|\dG_n^*|>0}\frac{\lambda_{\textrm{min}}(\bs{\Sigma}_{n})}{|\mathbb{T}_{n}^*|}
=\lambda_{\textrm{min}}(\bs{\Sigma})\indnex\hspace{1cm} \text{a.s.}
\end{equation*}
where $\lambda_{\textrm{min}}(\bs{A})$ denotes the smallest eigenvalue of matrix $\bs{A}$. Since $\bs{L}$ as well as $\bs{\Sigma}$ is definite positive, one has $\lambda_{\textrm{min}}(\bs{\Sigma})>0$. Therefore, as
\begin{equation*}
\|\widehat{\bs{\theta}}_{n}-\bs{\theta}\|^{2}\leq \frac{\cV_n}{\lambda_{\textrm{min}}(\bs{\Sigma}_{n-1})},
\vspace{-1ex}
\end{equation*}
we use Result~(\ref{th Mn1}) of Theorem~\ref{th Mn} to conclude that
\begin{equation*}
\ind{|\dG_n^*|>0}\|\widehat{\bs{\theta}}_{n}-\bs{\theta}\|^{2}= \cO \left(\frac{n}{|\dT_{n-1}^*|} \right)\indnex=
\cO \left(\frac{\log |\dT_{n-1}^*|}{|\dT_{n-1}^*|} \right)\indnex
\hspace{1cm}\text{a.s.}
\end{equation*}
which completes the proof of results~(\ref{thmaptheta1}).
\hspace{\stretch{1}}$ \Box$\\

We now prove the quadratic strong law (QSL).\\

\noindent\textbf{Proof of result~(\ref{thmaptheta2}) of Theorem~\ref{thmaptheta}:} The QSL is a direct consequence of result~(\ref{th Mn3}) of Theorem~\ref{th Mn} together with the fact that
$\wh{\bs{\theta}}_n-\bs{\theta}=\bs{\Sigma}_{n-1}^{-1}\bs{M}_n$. Indeed, we have
\begin{eqnarray*}
\lefteqn{\ind{|\dG_n^*|>0}\frac{1}{n}\sum_{\ell=1}^n\bs{M}_{\ell}^{t}\bs{\Sigma}_{\ell-1}^{-1} \bs{M}_{\ell}}\\
&=&\ind{|\dG_n^*|>0}\frac{1}{n}\sum_{\ell=1}^n (\wh{\bs{\theta}}_{\ell}-\bs{\theta})^{t}\bs{\Sigma}_{\ell-1}(\wh{\bs{\theta}}_{\ell}-\bs{\theta})\\
&=&\ind{|\dG_n^*|>0}\frac{1}{n}\sum_{\ell=1}^n |\dT_{\ell-1}^*|(\wh{\bs{\theta}}_{\ell}-\bs{\theta})^{t}
\ind{|\dG_{\ell-1}^*|>0}\frac{\bs{\Sigma}_{\ell-1}}{|\dT_{\ell-1}^*|}(\wh{\bs{\theta}}_{\ell}-\bs{\theta})\\
&=& \ind{|\dG_n^*|>0}\frac{1}{n}\sum_{\ell=1}^n |\dT_{\ell-1}^*|(\wh{\bs{\theta}}_{\ell}-\bs{\theta})^{t}
\bs{\Sigma}(\wh{\bs{\theta}}_{\ell}-\bs{\theta}) + o(1)
\hspace{1cm}\text{a.s.}
\end{eqnarray*}
which completes the proof.
\hspace{\stretch{1}}$ \Box$

\subsection{Strong consistency for the variance estimators}
\label{section main proof2}
For $n\geq1$, set
\begin{equation*}
\bs{V}_k=\left(\delta_{2k}\veps_{2k},\delta_{2k+1}\veps_{2k+1}\right)^t,\qquad
\wh{\bs{V}}_k=\left(\delta_{2k}\wh{\veps}_{2k},\delta_{2k+1}\wh{\veps}_{2k+1}\right)^t.
\end{equation*}
The almost sure convergence of $\wh{\sigma}^2_n$ and $\wh{\rho}_n$ is strongly related to that of $\wh{\bs{V}}_k-\bs{V}_k$.\\

\noindent\textbf{Proof of result~(\ref{apsigma1}) of Theorem~\ref{thmapsigmarho}:} Equation(\ref{apsigma1}) can be rewritten as
\begin{equation*}
\lim_{n\rightarrow\infty}\ind{|\dG_n^*|>0}\frac{1}{n}\sum_{k\in\mathbb{T}_{n-1}}\|\wh{\bs{V}}_k-\bs{V}_k\|^2
={4}(\pi-1)\sigma^2\indnex\hspace{1cm}\text{a.s.}
\end{equation*}
Once again, we are searching for a link between the sum of $\|\wh{\bs{V}}_k-\bs{V}_k\|$ and the processes $(\cA_n)$ and $(\cV_n)$ whose convergence properties were previously investigated. For $i\in\{0,1\}$ and $n\geq0$, let
\begin{equation*}
\bs{\Phi}_n^i= \left(\begin{array}{cccc}
\delta_{2(2^n)+i} & \delta_{2(2^n + 1)+i} & \cdots & \delta_{2(2^{n+1}-1)+i} \vspace{1ex}\\
\delta_{2(2^n)+i}X_{2^n} & \delta_{2(2^n + 1)+i}X_{2^n + 1} & \cdots & \delta_{2(2^{n+1}-1)+i}X_{2^{n+1}-1}
\end{array}\right)
\end{equation*}
be the collection of $(\delta_{2k+i}, \delta_{2k+i}X_k)^t$, $k\in\dG_n$, and set
\begin{equation*}
\bs{\Psi}_n= \left(\begin{array}{cc}
\bs{\Phi}_n^0&0\\
0&\bs{\Phi}_n^1
\end{array}\right).
\end{equation*}
Note that $\bs{\Psi}_n$ is a $4\times 2^{n+1}$ matrix. For all $n\geq 1$, we thus have, in the matrix form
\begin{eqnarray*}
\sum_{k\in\mathbb{G}_n}\|\wh{\bs{V}}_k-\bs{V}_k\|^2&=&\sum_{k\in\mathbb{G}_n}\delta_{2k}(\wh{\veps}_{2k}-\veps_{2k})^2+\delta_{2k+1}(\wh{\veps}_{2k+1}-\veps_{2k+1})^2,\\
&=&(\wh{\bs{\theta}}_n-\bs{\theta})^t\bs{\Psi}_n\bs{\Psi}_n^t(\wh{\bs{\theta}}_n-\bs{\theta}),\\
&=&\bs{M}_n^t\bs{\Sigma}_{n-1}^{-1}\bs{\Psi}_n\bs{\Psi}_n^t\bs{\Sigma}_{n-1}^{-1}\bs{M}_n,\\
&=&\bs{M}_n^t\bs{\Sigma}_{n-1}^{-1/2}\bs{\Delta}_n\bs{\Sigma}_{n-1}^{-1/2}\bs{M}_n,
\end{eqnarray*}
where
\begin{equation*}
\bs{\Delta}_n=\bs{\Sigma}_{n-1}^{-1/2}\bs{\Psi}_n\bs{\Psi}_n^t\bs{\Sigma}_{n-1}^{-1/2}=
\bs{\Sigma}_{n-1}^{-1/2}(\bs{\Sigma}_n-\bs{\Sigma}_{n-1})\bs{\Sigma}_{n-1}^{-1/2}.
\end{equation*}
Now, we can deduce from Corollary~(\ref{cv Gamman}) that
\begin{equation*}
\lim_{n\rightarrow\infty}\ind{|\dG_n^*|>0}\bs{\Delta}_n=(\pi-1)\rI_{4}\indnex\hspace{1cm}\text{a.s.}
\end{equation*}
which implies that
\begin{equation*}
\ind{|\dG_n^*|>0}\sum_{k\in\mathbb{G}_n}\|\wh{\bs{V}}_k-\bs{V}_k\|^2=\bs{M}_n^t\bs{\Sigma}_{n-1}^{-1}\bs{M}_n\left(\pi-1+o(1)\right)\ind{|\dG_n^*|>0}\hspace{0.2cm}\text{a.s.}
\end{equation*}
Therefore, we can conclude via convergence (\ref{th Mn3}) that
\begin{eqnarray*}
\lefteqn{\lim_{n\rightarrow\infty}\ind{|\dG_n^*|>0}\frac{1}{n}\sum_{k\in\mathbb{T}_{n-1}}\|\wh{\bs{V}}_k-\bs{V}_k\|^2}\\
&=&\lim_{n\rightarrow\infty}\ind{|\dG_n^*|>0}\frac{1}{n(\pi-1)}\sum_{\ell=1}^n \bs{M}_{\ell}^{t}\bs{\Sigma}_{\ell-1}^{-1}\bs{M}_{\ell}\ =\ 4(\pi-1)\sigma^2\indnex\hspace{1cm}\text{a.s.}
\end{eqnarray*}
which completes the proof.
\hspace{\stretch{1}}$ \Box$\\

\noindent\textbf{Proof of result~(\ref{apsigma2}) of Theorem~\ref{thmapsigmarho}:} First of all, one has
\begin{eqnarray*}
\wh{\sigma}^2_n-{\sigma}^2_n&=&\frac{1}{|\mathbb{T}_{n}^*|}
\sum_{k\in\mathbb{T}_{n-1}}\big(\|\wh{\bs{V}}_k\|^2-\|{\bs{V}}_k\|^2\big),\\
&=&\frac{1}{|\mathbb{T}_{n}^*|}\sum_{k\in\mathbb{T}_{n-1}}\big(\|\wh{\bs{V}}_k-{\bs{V}}_k\|^2+2(\wh{\bs{V}}_k-{\bs{V}}_k)^t\bs{V}_k\big).
\end{eqnarray*}
Set
\begin{equation*}
P_n=\sum_{k\in\mathbb{T}_{n-1}}(\wh{\bs{V}}_k-{\bs{V}}_k)^t\bs{V}_k=
\sum_{\ell=1}^{n}\sum_{k\in\mathbb{G}_{\ell-1}}(\wh{\bs{V}}_k-\bs{V}_k)^t\bs{V}_k.
\end{equation*}
We clearly have
\begin{equation*}
\Delta P_{n+1}=P_{n+1}-P_n=\sum_{k \in \dG_{n}}(\wh{\bs{V}}_k-\bs{V}_k)^t\bs{V}_k.
\end{equation*}
One can observe that for all $k \in \mathbb{G}_{n}$, $\wh{\bs{V}}_k-{\bs{V}}_k$ is $\cF_n^{\cO}$-measurable.
Consequently, $(P_n)$ is a real martingale transform for the filtration $\dF^{\cO}$.
Hence, we can deduce from the strong law of large numbers for martingale
transforms given in Theorem 1.3.24 of \cite{Duflo97} together with (\ref{apsigma1}) that
\begin{equation*}
\ind{|\dG_n^*|>0}P_n=o\left(\sum_{k\in\dT_{n-1}}||\wh{\bs{V}}_k-\bs{V}_k)||^2\right)=o(n)\hspace{1cm}\text{a.s.}
\end{equation*}
It ensures once again via convergence (\ref{apsigma1}) that
\begin{eqnarray*}
{\lim_{n\rightarrow\infty}\ind{|\dG_n^*|>0}\frac{|\dT_{n}^*|}{n}
(\wh{\sigma}^2_n-{\sigma}^2_n)}
&=&\lim_{n\rightarrow\infty}\ind{|\dG_n^*|>0}\frac{1}{n}\sum_{k\in\mathbb{T}_{n-1}}\|\wh{\bs{V}}_k-\bs{V}_k\|^2\\
&=&
4(\pi-1)\sigma^2\indnex\hspace{1cm}\text{a.s.}
\end{eqnarray*}
which  completes the proof of result~(\ref{apsigma2}).
\hspace{\stretch{1}}$ \Box$\\

\noindent\textbf{Proof of results~(\ref{aprho1}) and (\ref{aprho2}) of Theorem~\ref{thmapsigmarho}:}
We now turn to the study of the covariance estimator $\wh{\rho}_n$. We have
\begin{eqnarray*}
\wh{\rho}_n-{\rho}_n&=&\frac{1}{|\dT_{n-1}^{*01}|}\sum_{k\in\mathbb{T}_{n-1}}\delta_{2k}\delta_{2k+1}(\wh{\veps}_{2k}\wh{\veps}_{2k+1}-\veps_{2k}\veps_{2k+1}),\\
&=&\frac{1}{|\dT_{n-1}^{*01}|}\sum_{k\in\mathbb{T}_{n-1}}\delta_{2k}(\wh{\veps}_{2k}-\veps_{2k})\delta_{2k+1}(\wh{\veps}_{2k+1}-\veps_{2k+1})
+\frac{1}{|\dT_{n-1}^{*01}|}Q_n,
\end{eqnarray*}
where
\begin{eqnarray*}
 Q_n&=&\sum_{k\in\mathbb{T}_{n-1}}
\delta_{2k}\delta_{2k+1}(\wh{\veps}_{2k}-\veps_{2k})\veps_{2k+1}+\delta_{2k}\delta_{2k+1}(\wh{\veps}_{2k+1}-\veps_{2k+1})\veps_{2k}\\
&=&\sum_{k\in\mathbb{T}_{n-1}}(\wh{\bs{V}}_k-\bs{V}_k)^t\rJ_2\bs{V}_k,
\end{eqnarray*}
with
\begin{equation*}
\rJ_2=\left(\begin{array}{cc}0&1\\1&0\end{array}\right).
\end{equation*}
The process $(Q_n)$ is a real martingale transform for the filtration $\dF^{\cO}$ satisfying
\begin{equation*}
Q_n=o\left(\sum_{k\in\dT_{n-1}}||\wh{\bs{V}}_k-\bs{V}_k||^2\right)=o(n)\hspace{1cm}\text{a.s.}
\end{equation*}
It now remains to prove that
\begin{eqnarray}
\lefteqn{\lim_{n\rightarrow\infty}\ind{|\dG_n^*|>0}\frac{1}{n}
\sum_{k\in\dT_{n-1}}\delta_{2k}\delta_{2k+1}(\wh{\varepsilon}_{2k}-\varepsilon_{2k})(\wh{\varepsilon}_{2k+1}-\varepsilon_{2k+1})}\nonumber\\
&=&\lim_{n\rightarrow\infty}\frac{R_n}{n}\ =\rho(\pi-1) tr\big((\bs{L}^1)^{-1}(\bs{L}^{0,1})^2(\bs{L}^0)^{-1}\big)\indnex\label{cvgVJchap}
\hspace{1cm}\text{a.s.}
\end{eqnarray}
where
\begin{equation*}
R_n=\sum_{\ell=1}^n\bs{M}_{\ell}^t\bs{\Sigma}_{\ell-1}^{-1}(\bs{J}_2\otimes\bs{\Phi}_{\ell}^{01}(\bs{\Phi}_{\ell}^{01})^t)\bs{\Sigma}_{\ell-1}^{-1}\bs{M}_{\ell},
\end{equation*}
where $\otimes$ denotes the Kronecker product of matrices, i.e.
\begin{equation*}
\bs{J}_2\otimes\bs{\Phi}_{\ell}^{01}(\bs{\Phi}_{\ell}^{01})^t=\left(
\begin{array}{cc}
0&\bs{\Phi}_{\ell}^{01}(\bs{\Phi}_{\ell}^{01})^t\\
\bs{\Phi}_{\ell}^{01}(\bs{\Phi}_{\ell}^{01})^t&0
\end{array}
\right),
\end{equation*}
and $\bs{\Phi}_{\ell}^{01}$ is defined similarly as $\bs{\Phi}_{\ell}^{0}$ and $\bs{\Phi}_{\ell}^{1}$ by the collection of $(\delta_{2k}\delta_{2k+1}$, $\delta_{2k}\delta_{2k+1}X_k)^t$, $k\in\dG_{\ell}$. As $\bs{\Phi}_{n}^{01}(\bs{\Phi}_{n}^{01})^t=\bs{S}_n^{01}-\bs{S}_{n-1}^{01}$, proposition~\ref{mainlemma} implies that
\begin{equation*}
\lim_{n\rightarrow\infty}\bs{\Sigma}_{n-1}^{-1/2}(\bs{J}_2\otimes\bs{\Phi}_{n}^{01}(\bs{\Phi}_{n}^{01})^t)\bs{\Sigma}_{n-1}^{-1/2}=(\pi-1)\bs{\Sigma}^{-1/2}\bs{J}_2\otimes\bs{L}^{01}\bs{\Sigma}^{-1/2}\qquad\textrm{a.s.}
\end{equation*}
so that the asymptotic behavior of $R_n/n$ boils down to that of
\begin{equation*}
\sum_{\ell=1}^n\bs{M}_{\ell}^t\bs{\Sigma}_{\ell-1}^{-1/2}(\bs{J}_2\otimes\bs{L}^{01})\bs{\Sigma}_{\ell-1}^{-1/2}\bs{M}_{\ell}.
\end{equation*}
A proof along the same lines as in Section~\ref{section proof M} finally yields the expected results, i.e.
\begin{eqnarray*}
\lim_{n\rightarrow\infty}\ind{|\dG_n^*|>0}\frac{R_n}{n}
&=& \rho\frac{\pi-1}{\pi}tr\big((\bs{L}^1)^{-1}(\bs{L}^{0,1})^2(\bs{L}^0)^{-1}\big)\indnex\hspace{1cm}\text{a.s.}
\end{eqnarray*}
which completes the proof of convergence (\ref{cvgVJchap}). We then obtain
\begin{equation*}
\lim_{n\rightarrow\infty}\ind{|\dG_n^*|>0}\frac{|\dT^*_{n}|}{n}
(\wh{\rho}_n-{\rho}_n)=\rho\frac{\pi-1}{\bar{p}(1,1)}tr\big((\bs{L}^1)^{-1}(\bs{L}^{0,1})^2(\bs{L}^0)^{-1}\big)\indnex\hspace{0.5cm}\text{a.s.}
\end{equation*}
which completes the proof of Theorem~\ref{thmapsigmarho}.\hspace{\stretch{1}}$
\Box$

\subsection{Asymptotic normality}
\label{section th3}
Contrary to the previous literature on BAR processes, we cannot use the central limit theorem given by Propositions~7.8 and 7.9 of \cite{Ham94} as in \cite{Guy07, BSG09} because the normalizing term is now the number of observations and is therefore random. The approach used in \cite{DM08} strongly relies on the gaussian assumption for the noise sequence that does not hold here.
Instead, we use the central limit theorem for martingales
 given in Theorem~3.II.10 of Duflo \cite{Duflo97}. However, unlike the previous sections, this theorem can not be directly applied to the martingale $(\bs{M}_n)$ because the number of observed data in a given generation grows exponentially fast and the Lindeberg condition does not hold.  The solution is to use a new filtration. Namely, instead of using the observed generation-wise filtration, we will use the sister pair-wise one. Let
\begin{equation*}
\mathcal{G}^{\cO}_p=\cO \vee \sigma\{\delta_1X_1,\ (\delta_{2k}X_{2k}, \delta_{2k+1}X_{2k+1}),\ 1\leq k\leq p\}
\end{equation*}
be the $\sigma$-algebra generated by the whole history $\cO$ of the Galton-Watson process and all observed individuals up to the offspring of individual $p$. Hence $(\delta_{2k}\veps_{2k}, \delta_{2k+1}\veps_{2k+1})$ is $\mathcal{G}^{\cO}_k$-measurable.
In addition, assumptions {(\textbf{HN.1})} and  {(\textbf{HI})} imply that the processes $(\delta_{2k}\veps_{2k}, X_k\delta_{2k}\veps_{2k}, \delta_{2k+1}\veps_{2k+1}, X_k\delta_{2k+1}\veps_{2k+1})^t$, $(\delta_{2k}\veps_{2k}^2+\delta_{2k+1}\veps_{2k+1}^2-(\delta_{2k}+\delta_{2k+1})\sigma^2)$ and $(\delta_{2k}\delta_{2k+1}(\veps_{2k}\veps_{2k+1}-\rho))$ are $\mathcal{G}^{\cO}_k$-martingale difference sequences. In all the sequel, we will work under the probability $\mathbb{P}_{\overline{\cE}}$ and we denote by $\mathbb{E}_{\overline{\cE}}$ the corresponding expectation.\\

\noindent\textbf{Proof of Theorem~\ref{thmCLT}, first step: }We apply Theorem~3.II.10 of \cite{Duflo97}
to the $\mathcal{G}^{\cO}_k$-martingale $\bs{M}^{(n)}= (\bs{M}^{(n)}_p)_{\{p \geq 1\}} $ defined by 
\begin{equation*}
\bs{M}^{(n)}_p = \frac{1}{\sqrt{|\dT^*_n|}} \sum_{k=1}^{p}\bs{D}_k
\qquad \text{with} \qquad 
\bs{D}_k=\left( \begin{array}{cccc}
\delta_{2k}\veps_{2k}  \\
X_k\delta_{2k}\veps_{2k} \\
\delta_{2k+1}\veps_{2k+1} \\
X_k\delta_{2k+1}\veps_{2k+1}
\end{array}\right).
\end{equation*}
Set $\nu_n=|\dT_n|=2^{n+1}-1$. Note that if $k \notin \dT^*_n$, then $\bs{D}_k=0$ which implies that
\begin{equation*}
\bs{M}^{(n)}_{\nu_n}=  \frac{1}{\sqrt{|\dT^*_n|}} \sum_{k=1}^{|\dT_n|}\bs{D}_k= \frac{1}{\sqrt{|\dT^*_n|}} \sum_{k \in \dT^*_n}\bs{D}_k.
\end{equation*}
As the non-extinction set $\overline{\cE}$ is in $\mathcal{G}^{\cO}_{k}$ for every $k \geq 1$, it is easy to prove that 
\begin{eqnarray*}
\lefteqn{\mathbb{E}_{\overline{\cE}}[\bs{D}_k\bs{D}_k^t|\mathcal{G}^{\cO}_{k-1}]=\mathbb{E}[\bs{D}_k\bs{D}_k^t|\mathcal{G}^{\cO}_{k-1}]} \\
& = & \left(\begin{array}{cccc}
\sigma^2\delta_{2k} &   \sigma^2\delta_{2k}X_k & \rho  \delta_{2k}\delta_{2k+1} &  \rho  \delta_{2k}\delta_{2k+1}X_k  \\
 \sigma^2\delta_{2k}X_k &   \sigma^2\delta_{2k}X_k^2 & \rho  \delta_{2k}\delta_{2k+1}X_k &  \rho  \delta_{2k}\delta_{2k+1}X_k^2\\
\rho  \delta_{2k}\delta_{2k+1} &  \rho  \delta_{2k}\delta_{2k+1}X_k &  \sigma^2\delta_{2k+1} &   \sigma^2\delta_{2k+1}X_k\\
\rho  \delta_{2k}\delta_{2k+1}X_k &  \rho  \delta_{2k}\delta_{2k+1}X_k^2 &   \sigma^2\delta_{2k+1}X_k &   \sigma^2\delta_{2k+1}X_k^2
\end{array}\right),
\end{eqnarray*} 
and Corollary \ref{cv Gamman} gives  the $\mathbb{P}_{\overline{\cE}}$ almost sure limit of the increasing process
\begin{equation} \label{cvbr}
<\bs{M}^{(n)}>_{\nu_n}=  \frac{1}{|\dT^*_n|} \sum_{k \in \dT^*_n}\mathbb{E}_{\overline{\cE}}[\bs{D}_k\bs{D}^t_k|\mathcal{G}^{\cO}_{k-1}]  = \frac{\bs{\Gamma}_n }{|\dT^*_n|}  \xrightarrow[n\rightarrow\infty]{} \bs{\Gamma}.
\end{equation}
Therefore, the first assumption of Theorem~3.II.10 of \cite{Duflo97} holds under
$\mathbb{P}_{\overline{\cE}}$.  We now want to prove the Lindeberg condition that is the convergence  in probability to $0$ of the following expression $L_n$ for all $\epsilon >0$:
\begin{eqnarray*}
L_n & = & \frac{1}{|\dT^*_n|} \sum_{k \in \dT^*_n}\mathbb{E}_{\overline{\cE}}[\|\bs{D}_k\|^2 \ind{ \|\bs{D}_k\| > \epsilon \sqrt{|\dT^*_n|}} | \mathcal{G}^{\cO}_{k-1}] \\
& \leq & \frac{1}{|\dT^*_n|} \sum_{k \in \dT^*_n}\mathbb{E}_{\overline{\cE}}[\|\bs{D}_k\|^r| \mathcal{G}^{\cO}_{k-1}]  \dP_{\overline{\cE}}(  \|\bs{D}_k\| > \epsilon \sqrt{|\dT^*_n|} ~| \mathcal{G}^{\cO}_{k-1})\\
& \leq & \frac{\sup_{k\geq 0}
\dE[\|\bs{D}_k\|^r|\mathcal{G}^{\cO}_{k-1}]}{|\dT^*_n|} \sum_{k \in \dT^*_n} \frac{\mathbb{E}_{\overline{\cE}}[\|\bs{D}_k\|^2 ~| \mathcal{G}^{\cO}_{k-1}]}{{\epsilon^2 |\dT^*_n|}}
\end{eqnarray*}
{for some }$r>2$ and thanks to H\"{o}lder and Chebyshev inequalities. Besides, using Eq. (\ref{receq}) and similar calculations as in Lemma~\ref{lem SumX}, one readily obtains
\begin{equation*}
X_n^8\leq 2^7(1-\beta)^{-7}\sum_{k=0}^{r_n-1}\beta^k|\eta^8_{[\frac{n}{2k}]}|+2^7\beta^{8r_n}X_1^8.
\end{equation*}
Now, assumption (\textbf{HN.1}) together with $\beta<1$ yield the existence of a constant $C$ such that
\begin{equation*}
\sup_{k\geq 0}\mathbb{E}[X_k^8]\leq C(1+\mathbb{E}[X_1^8]),
\end{equation*}
and recall that $\mathbb{E}[X_1^8] < \infty$. Finally, since the entries of $\bs{D}_k$ are combinations of $\veps_{2k+i}$ and $X_k$, using again {(\textbf{HN.1})} and  {(\textbf{HI})}, one obtains that 
\begin{equation*}
\sup_{k\geq 0}
\dE[\|\bs{D}_k\|^r|\mathcal{G}^{\cO}_{k-1}]<\infty
\hspace{1cm}\text{a.s.}
\end{equation*} 
with $r=8$. The Lindeberg condition is thus proved, plugging the convergence (\ref{cvbr}) into the following equality:
\begin{equation*}
\frac{1}{|\dT^*_n|}\sum_{k \in \dT^*_n} \mathbb{E}_{\overline{\cE}}[\|\bs{D}_k\|^2 ~| \mathcal{G}^{\cO}_{k-1}] = tr\Big(\frac{1}{|\dT^*_n|} \sum_{k \in \dT^*_n} \mathbb{E}_{\overline{\cE}}[\bs{D}_k \bs{D}^t_k ~| \mathcal{G}^{\cO}_{k-1}]\Big) \xrightarrow[n\rightarrow\infty]{} tr(\bs{\Gamma}).
\end{equation*}
We can now conclude that under
$\mathbb{P}_{\overline{\cE}}$ 
\begin{equation*}
\frac{1}{\sqrt{|\dT^*_{n-1}|}}\sum_{k \in \mathbb{T}^*_{n-1}}\bs{D}_k=\frac{1}{\sqrt{|\dT^*_{n-1}|}}\bs{M}_n
\liml
\cN(0,\bs{\Gamma}).
\end{equation*}
Finally, result (\ref{CLTtheta}) follows from
Eq.~(\ref{thetadiff}) and Corollary~\ref{cv Gamman} together with Slutsky's Lemma.
\hspace{\stretch{1}}$\Box$\\

\noindent\textbf{Proof of Theorem~\ref{thmCLT}, second step:} On the one hand, we apply Theorem~3.II.10 of \cite{Duflo97}
to the $\mathcal{G}^{\cO}_p$-martingale $M^{(n)}= (M^{(n)}_p)_{\{p \geq 1\}} $ defined by 
\begin{equation*}
M^{(n)}_p = \frac{1}{\sqrt{|\dT^*_n|}} \sum_{k=1}^{p}v_k
\quad \text{and} \quad v_k=\delta_{2k}\veps_{2k}^2+\delta_{2k+1}\veps_{k+1}^2-(\delta_{2k}+\delta_{2k+1})\sigma^2.
\end{equation*}
As above, one clearly has
\begin{equation*}
M^{(n)}_{\nu_n}=  \frac{1}{\sqrt{|\dT^*_n|}} \sum_{k \in \dT^*_{n-1}}v_k= {\sqrt{|\dT^*_n|}} (\sigma_n^2-\sigma^2).
\end{equation*}
Using assumptions {(\textbf{HN.1})}, {(\textbf{HI})} and Lemma \ref{lemSumdelta} we compute the limit of the increasing process under $\mathbb{P}_{\overline{\cE}}$
\begin{equation*}
\lim_{n\rightarrow\infty}<M^{(n)}>_{\nu_n}=  (\tau^4-\sigma^4)+\frac{2\bar{p}(1,1)}{\pi}(\nu^2 \tau^4-\sigma^4) \qquad \mathbb{P}_{\overline{\cE}}\text{ a.s.}
\end{equation*}
Therefore, the first assumption of Theorem~3.II.10 of \cite{Duflo97} holds under
$\mathbb{P}_{\overline{\cE}}$.  Thanks to  assumptions {(\textbf{HN.2})} and  {(\textbf{HI})} we can prove  that for some $r>2$,
\begin{equation*}
\sup_{k\geq 0}
\dE_{\overline{\cE}}[\|v_k\|^r|\mathcal{G}^{\cO}_{k-1}]<\infty
\hspace{1cm}\text{a.s.}
\end{equation*} which implies the Lindeberg condition. 
Therefore, we obtain that under $\mathbb{P}_{\overline{\cE}}$
\begin{equation*}
\sqrt{|\dT^*_{n}|}(\sigma_n^2-\sigma^2)
\liml
\cN(0, (\tau^4-\sigma^4)+\frac{2\bar{p}(1,1)}{\pi}(\nu^2 \tau^4-\sigma^4) ).
\end{equation*}
Furthermore, we infer from Eq.~(\ref{apsigma2}) that
\begin{equation*}
\lim_{n\rightarrow\infty}\sqrt{|\dT^*_{n}|}(\wh{\sigma}_n^2-\sigma_n^2)=0\hspace{1cm}\mathbb{P}_{\overline{\cE}}\quad\text{a.s.}
\end{equation*}
which yields result~(\ref{CLTsigma}). \\

We turn now to the proof of result~(\ref{CLTrho}) with another $\mathcal{G}^{\cO}_p$-martingale $(M^{(n)})$ defined by
\begin{equation*}
M^{(n)}_p = \frac{1}{\sqrt{|\dT_{n-1}^{*01}|}} \sum_{k=1}^{p}\delta_{2k}\delta_{2k+1}(\veps_{2k}\veps_{2k+1}-\rho).
\end{equation*} 
As above, one easily shows that
\begin{equation*}
M^{(n)}_{\nu_n}=  \frac{1}{\sqrt{|\dT_{n-1}^{*01}|}} \sum_{i \in \dT^*_{n-1}}\delta_{2i}\delta_{2i+1}(\veps_{2i}\veps_{2i+1}-\rho)= {\sqrt{|\dT_{n-1}^{*01}|}} (\rho_n-\rho).
\end{equation*}
Using assumptions {(\textbf{HN.1})} and {(\textbf{H.I})}, we compute the limit of the increasing process
\begin{equation*}
\lim_{n\rightarrow\infty}<M^{(n)}>_{\nu_n}= \nu^2 \tau^4-\rho^2 \qquad \mathbb{P}_{\overline{\cE}}\text{ a.s.}
\end{equation*}
We also derive the Lindeberg condition. Consequently, we obtain that under $\mathbb{P}_{\overline{\cE}}$, one has
\begin{equation*}
\sqrt{|\dT^{*01}_{n-1}|}(\rho_n-\rho)
\liml
\cN(0,\nu^2 \tau^4-\rho^2).
\end{equation*}
Furthermore, we infer from (\ref{aprho2}) that
\begin{equation*}
\lim_{n\rightarrow\infty}\sqrt{|\dT^{*01}_{n-1}|}(\wh{\rho}_n-\rho_n)=0\hspace{1cm}\mathbb{P}_{\overline{\cE}}\quad\text{a.s.}
\end{equation*}
Finally, result~(\ref{CLTrho}) follows, which completes the proof of Theorem~\ref{thmCLT}.
\hspace{\stretch{1}}$\Box$

\appendix
\section{Quadratic strong law}
\label{appendixB}
In order to establish the quadratic strong law for the main martingale $(\bs{M}_n)$, we are going to study separately the asymptotic behavior of
$(\cW_n)$ and $(\cB_n)$ which appear in the main decomposition given by Equation~(\ref{maindecomart}).
\begin{Lemma}\label{lemlimW}
Under assumptions \emph{(\textbf{HN.1})}, \emph{(\textbf{HN.2})}, \emph{(\textbf{HO})} and \emph{(\textbf{HI})},  we have
\begin{equation*}
\lim_{n \rightarrow + \infty}\ind{|\dG_n^*|>0}\frac{1}{n} \cW_n= \frac{4(\pi-1)}{\pi} \sigma^2\indnex
\hspace{1cm}\text{a.s.}
\end{equation*}
\end{Lemma}

\noindent\textbf{Proof :} First of all, we have the decomposition $\cW_{n+1}=\cT_{n+1}+\cR_{n+1}$ where
\begin{eqnarray*}
\cT_{n+1}&=&\sum_{\ell=1}^n \frac{\Delta \bs{M}_{\ell+1}^t\bs{\Sigma}^{-1}\Delta \bs{M}_{\ell+1}}{|\dT_{\ell}^*|},\\
\cR_{n+1}&=&\sum_{\ell=1}^n \frac{\Delta \bs{M}_{\ell+1}^t(|\dT_{\ell}^*|\bs{\Sigma}_{\ell}^{-1}-\bs{\Sigma}^{-1})\Delta \bs{M}_{\ell+1}
}{|\dT_{\ell}^*|}.
\end{eqnarray*}
We first prove that
\begin{equation}\label{cv Tn}
\lim_{n \rightarrow + \infty}\ind{|\dG_n^*|>0}\frac{1}{n} \cT_n= \frac{4(\pi-1)}{\pi} \sigma^2\indnex
\hspace{1cm}\text{a.s.}
\end{equation}
As $\cT_n$ is a scalar and the trace is commutative, one can rewrite $\cT_{n+1}$ as $\cT_{n+1}=tr(\bs{\Sigma}^{-1/2}\bs{\cH}_{n+1}\bs{\Sigma}^{-1/2})$ where
\begin{equation*}
\bs{\cH}_{n+1}=\sum_{\ell=1}^{n} \frac{\Delta \bs{M}_{\ell+1}\Delta \bs{M}_{\ell+1}^t}{|\dT_{\ell}^*|}.
\end{equation*}
Our goal is to make use of the strong law of large numbers for martingale transforms, so we start by adding and
subtracting a term involving the conditional expectation of $\Delta\bs{\cH}_{n+1}$ given $\cF_n^{\cO}$. We have already seen in Section~\ref{section prop mart} that for all $n$, $\dE[\Delta \bs{M}_{n+1}\Delta \bs{M}_{n+1}^t|\cF_n^{\cO}]=\bs{\Gamma}_n-\bs{\Gamma}_{n-1}$. Consequently, we can split $\bs{\cH}_{n+1}$ into
two terms
\begin{equation*}
\bs{\cH}_{n+1}=\sum_{\ell=1}^{n} \frac{\bs{\Gamma}_{\ell}-\bs{\Gamma}_{\ell-1}}{|\dT_{\ell}^*|}+\bs{\cK}_{n+1},
\end{equation*}
where
\begin{equation*}
 \bs{\cK}_{n+1}=\sum_{\ell=1}^{n} \frac{\Delta \bs{M}_{\ell+1}\Delta \bs{M}_{\ell+1}^t-(\bs{\Gamma}_{\ell}-\bs{\Gamma}_{\ell-1})}{|\dT_{\ell}^*|}
\end{equation*}
On the one hand, it follows from Corollary~\ref{cv Gamman} and Lemma~\ref{lemcard} that
\begin{equation*}
\lim_{n \rightarrow + \infty}\ind{|\dG_{n}^*|>0}\frac{\bs{\Gamma}_n-\bs{\Gamma}_{n-1}}{|\dT_n^*|} = \frac{\pi-1}{\pi}\bs{\Gamma}\indnex\hspace{1cm}\text{a.s.}
\end{equation*}
Thus, Cesaro convergence and the remark that $\{|\dG_{\ell}^*|=0\} \subset \{|\dG_{n}^*|=0\}$ for all $\ell \le n$ yield
\begin{eqnarray*}
\lim_{n \rightarrow + \infty}\ind{|\dG_{n}^*|>0}\frac{1}{n}\sum_{\ell=1}^{n} \frac{\bs{\Gamma}_{\ell}-\bs{\Gamma}_{\ell-1}}{|\dT_{\ell}^*|}
&=& \lim_{n \rightarrow + \infty}\ind{|\dG_{n}^*|>0}\frac{1}{n}\sum_{\ell=1}^{n} \ind{|\dG_{\ell}^*|>0}\frac{\bs{\Gamma}_{\ell}-\bs{\Gamma}_{\ell-1}}{|\dT_{\ell}^*|}\\
&=&\frac{\pi-1}{\pi}\bs{\Gamma}\indnex\hspace{1cm}\text{a.s.} 
\end{eqnarray*}
On the other hand, the sequence $(\bs{\cK}_n)$ is obviously a matrix martingale transform and tedious but straightforward calculations, together with Lemmas~\ref{lem SumX} and  \ref{lemSumX4} and the strong law of large numbers for martingale transforms given in Theorem 1.3.24 of \cite{Duflo97} imply that
$\ind{|\dG_{n}^*|>0}\bs{\cK}_n=o(n)$ a.s.
Hence, we infer from the equation above that
\begin{equation} \label{cv Hn}
\lim_{n \rightarrow + \infty}\ind{|\dG_{n}^*|>0}\frac{1}{n}\bs{\cH}_n= \frac{\pi-1}{\pi}\bs{\Gamma}\indnex
\hspace{1cm}\text{a.s.}
\end{equation}
Finally, we obtain
\begin{equation*}
\lim_{n \rightarrow + \infty}\ind{|\dG_{n}^*|>0}\frac{1}{n}\cT_n
= \frac{\pi-1}{\pi}tr(\bs{\Sigma}^{-1/2}\bs{\Gamma}\bs{\Sigma}^{-1/2})\indnex
 = \frac{\pi-1}{\pi}4\sigma^2\indnex
\hspace{.3cm}\text{a.s.}\\
\end{equation*}
which proves (\ref{cv Tn}).
We now turn to the asymptotic behavior of $\cR_{n+1}$. We know from Proposition~\ref{mainlemma} that $\ind{|\dG_n^*|>0}(|\dT_n^*|\bs{\Sigma}_n^{-1}-\bs{\Sigma}^{-1})$ goes to $0$ as $n$ goes to infinity. Hence, for all positive $\epsilon$ and for large enough $n$, one has
\begin{equation*}
\ind{|\dG_n^*|>0} |\Delta \bs{M}_{n+1}^t(|\dT_n^*|\bs{\Sigma}_n^{-1}-\bs{\Sigma}^{-1})\Delta \bs{M}_{n+1}|
\leq \ind{|\dG_n^*|>0}4\epsilon \Delta \bs{M}_{n+1}^t\Delta \bs{M}_{n+1}.
\end{equation*}
Using again that $\{|\dG_{\ell}^*|=0\} \subset \{|\dG_{n+1}^*|=0\}$ for all $\ell \le n+1$, we rewrite 
\begin{equation*}
\ind{|\dG_{n+1}^*|>0}\cR_{n+1}\!=\!\ind{|\dG_{n+1}^*|>0}\!\!\sum_{\ell=1}^n\! \ind{|\dG_{\ell}^*|>0}\!\frac{\Delta \bs{M}_{\ell+1}^t(|\dT_{\ell}^*|\bs{\Sigma}_{\ell}^{-1}-\bs{\Sigma}^{-1})\Delta \bs{M}_{\ell+1}}{|\dT_{\ell}^*|}.
\end{equation*}
Hence,
\begin{eqnarray*}
\ind{|\dG_{n+1}^*|>0}|\cR_{n+1}| &\le& 4\epsilon \,\ind{|\dG_{n+1}^*|>0}\!\!\sum_{\ell=1}^n\! \ind{|\dG_{\ell}^*|>0}\!\frac{\Delta \bs{M}_{\ell+1}^t \Delta \bs{M}_{\ell+1}}{|\dT_{\ell}^*|}\\
 &\le& 4\epsilon \,\ind{|\dG_{n+1}^*|>0} \sum_{\ell=1}^n\! \frac{\Delta \bs{M}_{\ell+1}^t \Delta \bs{M}_{\ell+1}}{|\dT_{\ell}^*|} \\
&\le& 4\epsilon \,\ind{|\dG_{n+1}^*|>0}\, tr(\bs{\cH}_{n+1}).
\end{eqnarray*}
This last inequality holding for any positive  $\epsilon$ and large enough $n$, the limit given by Equation~(\ref{cv Hn}) entails that
\begin{equation*}
\lim_{n \rightarrow + \infty}\ind{|\dG_n^*|>0}\frac{1}{n} \cR_n= 0\hspace{1cm}\text{a.s.}
\end{equation*}
which completes the proof of Lemma~\ref{lemlimW}.
\hspace{\stretch{1}}$ \Box$

\begin{Lemma}
\label{lemlimB}
Under assumptions \emph{(\textbf{HN.1})}, \emph{(\textbf{HN.2})}, \emph{(\textbf{HO})} and \emph{(\textbf{HI})}, we have
\begin{equation*}
\lim_{n \rightarrow + \infty}\ind{|\dG_n^*|>0}\frac{1}{n} \cB_n= 0\hspace{1cm}\text{a.s.}
\end{equation*}
\end{Lemma}

\noindent\textbf{Proof :} The result is obvious on the extinction set $\cE$. Now let us work on $\overline{\cE}$. 
Now for $i\in\{0,1\}$ and $n\geq1$, let
$
\bs{\xi}_n^i= \left(\veps_{2^n+i},
\veps_{2^n+2+i},
\ldots,
\veps_{2^{n+1}-2+i}
\right)^t$,
be the collection of $\veps_k$, $k\in\dG_n^i$, and set
$\bs{\xi}_n= \left(
\bs{\xi}_n^0,
\bs{\xi}_n^1
\right)^t$.
Note that $\bs{\xi}_n$ is a column vector of size $2^{n+1}$. 
With these notation, one has
\begin{equation*}
\cB_{n+1}=2\sum_{\ell=1}^n \bs{M}_{\ell}^t\bs{\Sigma}_{\ell}^{-1}\Delta \bs{M}_{\ell+1}=2\sum_{\ell=1}^n \bs{M}_{\ell}^t\bs{\Sigma}_{\ell}^{-1}\bs{\Psi}_{\ell}\bs{\xi}_{\ell+1}.
\end{equation*}
The sequence $(\cB_n)$ is a real martingale transform satisfying
\begin{equation*}
\Delta \cB_{n+1}=\cB_{n+1}- \cB_{n}=2\bs{M}_n^t\bs{\Sigma}_n^{-1}\bs{\Psi}_n\bs{\xi}_{n+1}.
\end{equation*}
Consequently, via the strong law of large numbers for martingale transforms, we find that either $(\cB_n)$ converges a.s. or $\cB_{n}=o(<\cB>_n)$ a.s. where
\begin{equation*}
<\cB>_{n+1}= 4\sum_{\ell=1}^n \bs{M}_{\ell}^t\bs{\Sigma}_{\ell}^{-1}\bs{\Psi}_{\ell}\bs{C}\bs{\Psi}_{\ell}^t\bs{\Sigma}_{\ell}^{-1}\bs{M}_{\ell},
\end{equation*}
with
\begin{equation*}
\bs{C}=\left(\begin{array}{cc}
\sigma^2&\rho\\
\rho&\sigma^2
\end{array}\right)\otimes \rI_{2^n}.
\end{equation*}
As $\bs{C}$ is definite positive under assumption (\textbf{HN.1}), one has $\bs{C}\leq2\sigma^2\rI_{2^{n+1}}$ in the sense that $2\sigma^2\rI_{2^{n+1}}-\bs{C}$ is semi definite positive. Hence, one has
\begin{equation*}
<\cB>_{n+1}\leq 8\sigma^2\sum_{\ell=1}^n \bs{M}_{\ell}^t\bs{\Sigma}_{\ell}^{-1}\bs{\Psi}_{\ell}\bs{\Psi}_{\ell}^t\bs{\Sigma}_{\ell}^{-1}\bs{M}_{\ell}.
\end{equation*}
Now, by definition, one has
\begin{equation*}
\bs{\Sigma}_{\ell}^{-1}\bs{\Psi}_{\ell}\bs{\Psi}_{\ell}^t\bs{\Sigma}_{\ell}^{-1}=\left(
\begin{array}{cc}
(\bs{S}_{\ell}^0)^{-1}\bs{\Phi}_{\ell}^0(\bs{\Phi}_{\ell}^0)^t(\bs{S}_{\ell}^0)^{-1}&0\\
0&(\bs{S}_{\ell}^1)^{-1}\bs{\Phi}_{\ell}^1(\bs{\Phi}_{\ell}^1)^t(\bs{S}_{\ell}^1)^{-1}
\end{array}
\right).
\end{equation*}
We now use Lemma~B.1 of \cite{BSG09} on each entry to obtain
\begin{equation*}
\bs{\Sigma}_{\ell}^{-1}\bs{\Psi}_{\ell}\bs{\Psi}_{\ell}^t\bs{\Sigma}_{\ell}^{-1}\leq  \bs{\Sigma}_{\ell-1}^{-1}-\bs{\Sigma}_{\ell}^{-1},
\end{equation*}
as the matrix $l_k$ in that lemma is definite positive. Therefore, we obtain that
\begin{equation*}
<\cB>_{n+1} \leq 8\sigma^2\sum_{\ell=1}^n\bs{M}_{\ell}^t(\bs{\Sigma}_{\ell-1}^{-1}-\bs{\Sigma}_{\ell}^{-1})\bs{M}_{\ell}=8\sigma^2\cA_n.
\end{equation*}
Finally, we deduce from the main decomposition given by Equation~(\ref{maindecomart}) and Lemma~\ref{lemlimW} that
\begin{equation*}
\ind{|\dG_n^*|>0}(\cV_{n+1}+ \cA_n)=o(\cA_n)+\cO(n)
\hspace{1cm}\text{a.s.}
\end{equation*}
leading to $\ind{|\dG_n^*|>0}\cV_{n+1}=\cO(n)$ and $\ind{|\dG_n^*|>0}\cA_{n}=\cO(n)$ a.s. as $\cV_{n+1}$ and $\cA_{n}$ are non-negative. This implies in turn that $\ind{|\dG_n^*|>0}\cB_n=o(n)$ a.s. completing
the proof of Lemma~\ref{lemlimB}.\hspace{\stretch{1}}$ \Box$

\section{Wei's Lemma}
\label{appendixC}
\noindent In order to prove Proposition~\ref{lemlimM}, we shall apply Wei's Lemma given in \cite{Wei87} page 1672, to each entry of the vector-valued main martingale
\begin{equation*}
\bs{M}_n= \sum_{\ell=1}^n \sum_{k \in \dG_{\ell-1}} \left(
\delta_{2k}\veps_{2k},
\delta_{2k}X_k\veps_{2k},
\delta_{2k+1}\veps_{2k+1},
\delta_{2k+1}X_k\veps_{2k+1}
\right)^t.
\end{equation*}
For $i\in\{0,1\}$, denote
\begin{equation*}
P_n^i=\sum_{\ell=1}^n \sum_{k \in \dG_{\ell-1}}\delta_{2k+i}\veps_{2k+i}
\hspace{1cm}\text{and}\hspace{1cm}
Q_n^i=\sum_{\ell=1}^n \sum_{k \in \dG_{\ell-1}}\delta_{2k+i}X_k\veps_{2k+i}.
\end{equation*}
On the set $\overline{\cE}$, these processes can be rewritten as
\begin{equation*}
P_n^i=\sum_{\ell=1}^n\sqrt{|\dG_{\ell-1}^*|}v_{\ell}^i,\qquad Q_n^i=\sum_{\ell=1}^n\sqrt{|\dG_{\ell-1}|}w_{\ell}^i,
\end{equation*}
where
\begin{eqnarray*}
v_n^i&=&\ind{|\dG_{n-1}^*|>0}\frac{1}{\sqrt{|\dG_{n-1}^*|}}\sum_{k \in \dG_{n-1}}\delta_{2k+i}\veps_{2k+i},\\
w_n^i&=&\ind{|\dG_{n-1}^*|>0}\frac{1}{\sqrt{|\dG_{n-1}^*|}}\sum_{k \in \dG_{n-1}}\delta_{2k+i}X_k\veps_{2k+i}.
\end{eqnarray*}
On the one hand, we clearly have $\dE[v_{n+1}^i|\cF_n^{\cO}]=0$ and $\dE[(v_{n+1}^i)^2|\cF_n^{\cO}]=\sigma^2\frac{Z_{n+1}^i}{|\dG_{n}^*|}$ a.s. on $\overline{\cE}$.
Moreover, it follows from Cauchy-Schwarz inequality that
\begin{eqnarray*}
\dE[(v_{n+1}^i)^4|\cF_n^{\cO}]
&=&\frac{\ind{|\dG_{n}^*|>0}}{|\dG_{n}^*|^2}\!\!\sum_{k \in \dG_{n}}\delta_{2l+i}\dE[\veps_{2k+i}^4|\cF_n^{\cO}]
\\
&&+\frac{\ind{|\dG_{n}^*|>0}}{|\dG_{n}^*|^2}\!\!\sum_{p \in \dG_{n}}\sum_{k\neq p}\delta_{2p+i}\delta_{2k+i}\dE[\veps_{2p+i}^2|\cF_n^{\cO}]\dE[\veps_{2k+i}^2|\cF_n^{\cO}]\\
&\leq & 3C\ind{|\dG_{n}^*|>0}\sup_{k\in \dG_n}\dE[\veps_{2k+i}^4|\cF_n^{\cO}]
\hspace{0.5cm}\text{a.s.}
\end{eqnarray*}
as $Z^i_{n+1}|\dG_{n}^*|^{-1}$ is bounded.This implies that $\sup \dE[(v_{n+1}^i)^4|\cF_n^{\cO}]< +\infty$ a.s.
Consequently, we deduce from Wei's Lemma that for all $\eta>1/2$,
\begin{equation*}
\ind{|\dG_{n-1}^*|>0}(P_n^i)^2=o(|\dT_{n-1}^*| n^\eta)\indnex
\hspace{1cm}\text{a.s.}
\end{equation*}
On the other hand,
it is not hard to see that $\dE[w_{n+1}^i|\cF_n^{\cO}]=0$ a.s.
Moreover, it follows from Cauchy-Schwarz inequality that,
\begin{eqnarray*}
\lefteqn{\dE[(w_{n+1}^i)^4|\cF_n^{\cO}]}\\
&\leq&\frac{\ind{|\dG_{n}^*|>0}}{|\dG_{n}^*|^2}\left(\sum_{k \in \dG_{n}}
\delta_{2k+i}X_{k}^4\dE[\veps_{2k+i}^4|\cF_n^{\cO}]
+\sigma^4\sum_{p \in \dG_{n}}\sum_{k\neq p}\delta_{2p+i}\delta_{2k+i}X_p^2X_k^2\right)\\
&\!\!\leq \!\!& 3\ind{|\dG_{n}^*|>0}\left(\sup_{k\in \dG_n}\dE[\veps_{2k+i}^4|\cF_n^{\cO}]\right)
\left(\frac{1}{|\dG_{n}^*|}\!\!\sum_{k\in \dG_{n}}\delta_{2k+i}X_l^2\right) ^2
\hspace{0.5cm}\text{a.s.}
\end{eqnarray*}
which is finite from Proposition~\ref{propSumX2}.
We deduce from Wei's Lemma applied to $Q_n^i$ that for all $\eta>1/2$,
$\ind{|\dG_{n-1}^*|>0}\|Q_n^i\|^2=o(|\dT_{n-1}^*| n^\eta)
$
a.s. which completes the proof of Proposition~\ref{lemlimM}.\hspace{\stretch{1}}$ \Box$

\bibliographystyle{plain}
\bibliography{blabla}

\end{document}